\documentclass[oneside]{article}
\usepackage{amsfonts}
\usepackage[a4paper, total={7in, 9in}]{geometry}
\usepackage{bm}
\usepackage{amsmath}
\usepackage{hyperref}
\usepackage{multicol}
\usepackage{authblk}
\usepackage{amsthm}
\usepackage{graphicx}
\usepackage{float}
\usepackage[normalem]{ulem}

\usepackage{pgfplots}
\usepackage{subcaption}
\usepackage{pdflscape}
\usepackage{array}
\usepackage{epstopdf}
\usepackage{cleveref}
\usepackage{amsmath}
\usepackage{amsfonts}
\usepackage{amssymb}
\mathchardef\mhyphen="2D

\newcommand{\blue}{}

\usepackage{authblk}
\usepackage{amsthm}

\usepackage{tikz}
\usepackage{tikz-network}

\begin{document}

\title{Regularity-Conforming Neural Networks (ReCoNNs) for solving Partial Differential Equations}
\author[1]{Jamie M. Taylor (\texttt{jamie.taylor@cunef.edu}) }
\author[2,3,4]{David Pardo (\texttt{david.pardo@ehu.eus})}
\author[3,5]{Judit Mu\~noz-Matute (\texttt{jmunoz@bcamath.org})}
\affil[1]{CUNEF Universidad, Madrid, Spain }
\affil[2]{University of the Basque Country (UPV/EHU), Leioa, Spain}
\affil[3]{Basque Center for Applied Mathematics (BCAM), Bilbao, Spain}
\affil[4]{Ikerbasque: Basque Foundation for Science, Bilbao, Spain}
\affil[5]{Oden Institute for Computational Engineering and Sciences, The University of Texas at Austin, USA}
\date{}
\maketitle

\abstract{Whilst the Universal Approximation Theorem guarantees the existence of approximations to Sobolev functions --the natural function spaces for PDEs-- by Neural Networks (NNs) of sufficient size, low-regularity solutions may lead to poor approximations in practice. For example, classical fully-connected feed-forward NNs fail to approximate continuous functions whose gradient is discontinuous when employing strong formulations like in Physics Informed Neural Networks (PINNs). In this article, we propose the use of regularity-conforming neural networks, where {\it a priori} information on the regularity of solutions to PDEs can be employed to construct proper architectures. We illustrate the potential of such architectures via a two-dimensional (2D) transmission problem, where the solution may admit discontinuities in the gradient across interfaces, as well as power-like singularities at certain points. In particular, we formulate the weak transmission problem in a PINNs-like strong formulation with interface {\blue and continuity} conditions. Such architectures are partially explainable; discontinuities are explicitly described,  allowing the introduction of novel terms into the loss function. We demonstrate via several model problems in one and two dimensions the advantages of using regularity-conforming architectures in contrast to classical architectures. The ideas presented in this article {\blue easily extend} to problems in higher dimensions.\\

\noindent\textit{Keywords}: Regularity-conforming, neural networks, interface condition, transmission problem, singularities, PINNs.
}

\tableofcontents

\section{Introduction}
There has been a recent wealth of work on the use of Neural Networks (NNs) for numerically solving Partial Differential Equations (PDEs). Their advantages are particularly seen in problems in high dimensions \cite{han2018solving,zhang2022fbsde}, where quantitative versions of the Universal Approximation Theorem \cite{hornik1990universal} can guarantee good approximation properties without the curse of dimensionality. They are also capable to learn solutions to parametric problems \cite{aldirany2024operator,khoo2021solving,uriarte2022finite,brevis2024learning}, which is advantageous in inverse problems. 

In the last decade, many strategies have been developed to define loss functions to be minimized during training, such as the Physics Informed Neural Networks (PINNs) \cite{cai2021physics,cuomo2022scientific,de2024physics}, and all its variants \cite{pang2019fpinns,yang2021b,mahmoudabadbozchelou2022nn,jagtap2020extended}, the Variational PINNs \cite{kharazmi2019variational,kharazmi2021hp,berrone2022variational}, the Robust Variational PINNs \cite{rojas2024robust}, the Deep Ritz method \cite{yu2018deep}, the Deep Double Ritz method \cite{uriarte2023deep}, the Galerkin NNs \cite{ainsworth2021galerkin,ainsworth2022galerkin}, and the Deep Fourier method \cite{taylor2023deep,taylor2023deep2,TAYLOR2024116997}. All these techniques are based on the classical theory of variational methods for approximating PDEs. Similarly, several NN architectures have been considered in the literature, being feed-forward NNs with multiple hidden layers \cite{raissi2019physics} the most popular choice for approximating PDEs (see Figure \ref{figFCNN}). Other examples of architectures include convolutional NNs \cite{gao2021phygeonet} suitable for rectangular domains with uniform grids, multi-level networks that reduce optimisation errors \cite{aldirany2024multi}, recurrent NNs when temporal dynamics are involved \cite{ren2022phycrnet}, residual NNs \cite{cheng2021deep} to avoid gradient vanishing problems, autoencoders \cite{bhattacharya2021model} to reduce dimensionality, and graph NNs \cite{li2020multipole} for problems involving non-Euclidean structures. Finally, the activation functions \cite{sharma2017activation} such as ReLU or $tanh$ play a crucial role in NNs, introducing nonlinearity into the system, which allows NNs to learn complex dynamics.

Certain PDEs models {\blue exhibit} only {\blue limited} regularity, {\blue showing} singular behaviour over lower dimensional sets. Although Universal Approximation guarantees that sufficiently large NNs are capable of approximating such irregular solutions in appropriate Sobolev norms, this does not in general imply that the deep learning algorithm will lead to a good approximation or stable convergence towards a solution. In fact, there are known cases of problems in machine learning where, despite the existence of NNs with sufficient approximability and well-defined loss functions, no algorithm may succesfully obtain an approximation within a certain tolerance {\blue error} \cite{colbrook2022difficulty}. 

{\blue The general approach for solving a PDE with NNs is to define an architecture, consider a loss function that is minimised at the solution of the PDE, and employ an appropriate optimisation strategy to obtain the minimiser. The focus of this work is on the firs aspect: To define appropriate architectures for solving PDEs with NNs}. As we will illustrate in Section 2, classical feed-forward fully connected NNs often fail to successfully approximate solutions to PDEs with low regularity. One may be tempted to use smooth activation functions to approximate functions in $H^1$. But due to discontinuities in the gradient, this choice will lead to Gibb's-type phenomenon, which in turn produce large integration errors, leading to poor convergence properties. On the other hand, low-regularity activation functions, such as ReLU, lead to numerical instabilities when the loss function involves first-order derivatives of the NN itself \cite{magueresse2024adaptive}.

In this article, we propose a solution to {\blue solve low-regularity PDEs} by introducing {\it Regularity-Conforming NNs (ReCoNNs)}, where {\it a priori} knowledge on the regularity of the solutions to the PDE allows the design of proper architectures {\blue that lead to fast convergence and good-accuracy solutions}. To illustrate the idea of ReCoNNs, we focus on the particular case of the two-dimensional transmission problem with discontinuous materials. The choice of this model problem is two-fold: (a) precise regularity statements are known from the literature \cite{blum1982finite,brenner1999multigrid,brenner1997multigrid,brenner2003multigrid}, allowing us to easily define appropriate architectures , and (b) the problem is of practical interest, as it is an elliptic problem capable of modelling phenomena in domains that possess sharp material interfaces and are of great interest in several applications \cite{bonnafont2024finite,ong2005three,tyni2024boundary}. The extension of the ideas presented in this article to problems in three dimensions is straightforward.

To construct the proposed ReCoNNs, we use theoretical results on the location and type of singularity, which are available in multiple applications (e.g., the location often corresponds to the material discontinuities). We then train the ReCoNN architecture to efficiently approximate the specific form (exponents, stress intensity factors, etc.) of the singularity, which is often unknown {\it a priori}. In this way, we use only available information while we reconstruct the unknown part. In addition, ReCoNNs enjoy partial {\it explainability} of the obtained solution. That is, the network provides not only approximation of the singularities, but also describes them.  

In the scenarios defined in this article, we consider loss functions in strong form similar to the classical PINNs. For solutions that are smooth away from point singularites, a classical PINNs implementation is possible. However, selecting smooth feed-forward NNs generally leads to significant numerical instabilities or Gibbs phenomena near the singularity that lead to poor approximations. We show that the ReCoNNs architecture is able to well-approximate the solution in these cases. The partial explainability property of ReCoNNs allows us to consider a wider range of loss functions for our PDEs, where we may include, for example, interface conditions that are characteristic of weak formulations, whilst implementing the strong-form PDE. It is rather simple with ReCoNNs to extend the loss functional in strong form with explicit interface conditions to obtain a proper PINN-type loss to solve low-regularity problems that cannot be solved with current PINN methods \cite{kharazmi2021hp}. In addition, in practice we observe a good convegence towards the target function.

The structure of the paper is as follows. In Section \ref{secMotivation}, we first motivate the problem via a simple one-dimensional example by considering the approximation of a function in $H^1\setminus C^1$, which corresponds to the weak solution of an elliptic ODE. We demonstrate how low-regularity (ReLU) and high-regularity ($tanh$) activation functions in a fully-connected feed-forward NN {\blue can} fail to approximate effectively such a solution in the $H^1$-norm. We demonstrate how a ReCoNN architecture successfully approximates the target function, which corresponds to the solution of a weak-form PDE.

After presenting the key ideas in a simplified 1D setting, in Section \ref{secModel}, we turn to the 2D transmission problem with discontinuous materials, defining the problem and identifying the two types of singular behaviour encountered: jump discontinuities of the normal component of the gradient across interfaces and power-type singularities at singular vertices. The latter can be encountered in two scenarios: the case of a polygonal domain that admits re-entrant corners and when the material interfaces share a common vertex. We define appropriate ReCoNNs architectures for all these cases by adding singular components to the classical feed-forward NN in order to capture the singular behaviours of the solution. 

In the numerical results in Section \ref{secNumerics}, we test the introduced ReCoNNs architectures in each scenario. We first consider a case which only admits discontinuities in the gradient across interfaces. Here, we  introduce a loss function that allows a strong-form implementation with an interface condition in order to approximate the weak solution. We present numerical results showing strong approximation capabilities of the ReCoNNs architecture. We then consider the classical L-shaped domain problem, which admits a power-like singularity at its singular vertex. Furthermore, as the target function is $C^2$ on its domain, we compare our results to a PINNs formulation with a classical architecture, demonstrating the superior approximation capability of our proposed approach. Finally, we consider the transmission problem in a square domain {\blue with four quadratnts, each corresponding to} a distinct material. Its solution admits both types of singular behaviours simultaneously. Again, we define an appropriate ReCoNN architecture and implement a strong-form PDE with interface condition to approximate solutions, showing good approximation capability. As the interface condition impedes the use of a classical-PINNs type loss, we compare a classical architecture to the conforming architecture by employing the $H^1$-norm of the error, again, observing far stronger approximation capabilities {\blue with ReCoNNs.} 

Finally, Section \ref{secConclusions} summarizes the results of the article and the possible avenues for future investigation. Appendix \ref{App} explains the construction of the exact solution for an interior vertex problem we considered in the numerical results.

\section{Difficulty of NNs to approximate jumps in the gradient}\label{secMotivation}

\subsection{The architecture of a classical NN}

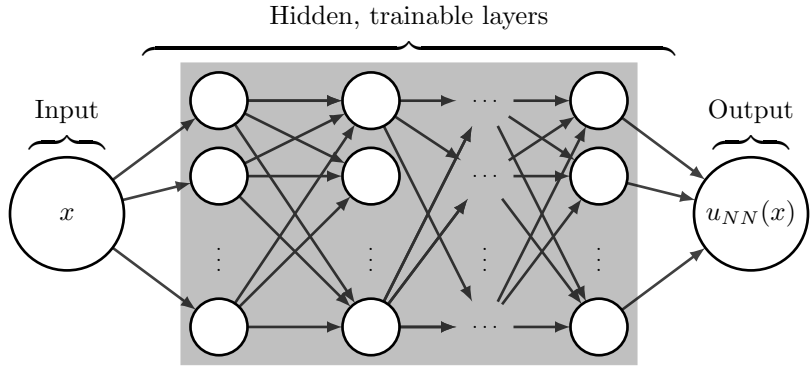
\begin{figure}[H]

\begin{center}
\begin{tikzpicture}

\filldraw [fill=black, draw=black,opacity=0.25] (1.5,-0.5) rectangle (7.5,3.5);

\draw (0,2.85) node{Input};
\draw (0,2.5) node {$\overbrace{\hspace{1cm}}^{\hspace{1cm}}$} ;

\draw (4.5,4.1) node {Hidden, trainable layers};
\draw (4.5,3.75) node {$\overbrace{\hspace{7cm}}^\text{\hspace{1cm}}$} ;

\draw (9,2.85) node {Output};
\draw (9,2.5) node {$\overbrace{\hspace{1cm}}^\text{\hspace{1cm}}$};

\Vertex[x=0,y=1.5,label=$x$,color=white,size=1.5,fontsize=\normalsize]{X}
\Vertex[x=2,y=3,color=white,size=0.75]{A11}
\Vertex[x=4,y=3,color=white,size=0.75]{A21}
\Vertex[x=7,y=3,color=white,size=0.75]{A31}
\Vertex[x=2,y=2,color=white,size=0.75]{A12}
\Vertex[x=4,y=2,color=white,size=0.75]{A22}
\Vertex[x=7,y=2,color=white,size=0.75]{A32}
\Vertex[x=2,y=0,color=white,size=0.75]{A13}
\Vertex[x=4,y=0,color=white,size=0.75]{A23}
\Vertex[x=7,y=0,color=white,size=0.75]{A33}

\Vertex[x=9,y=1.5,label=$u_{NN}(x)$,color=white,size=1.5,fontsize=\normalsize]{tildeu}

\Vertex[x=5.5,y=3,style={color=black,opacity=0.0},label=$\color{black}\hdots$,size=0.75]{Ah1}
\Vertex[x=5.5,y=2,style={color=white,opacity=0},label=$\color{black}\hdots$,size=0.75]{Ah2}
\Vertex[x=5.5,y=0,style={color=white,opacity=0},label=$\color{black}\hdots$,size=0.75]{Ah3}

\Edge[Direct,lw=1pt](X)(A11)
\Edge[Direct,lw=1pt](X)(A12)
\Edge[Direct,lw=1pt](X)(A13)

\Edge[color=white,label={$\color{black}\vdots$}](A12)(A13)
\Edge[color=white,label={$\color{black}\vdots$}](A22)(A23)
\Edge[color=white,label={$\color{black}\vdots$}](Ah2)(Ah3)
\Edge[color=white,label={$\color{black}\vdots$}](A32)(A33)

\Edge[Direct,lw=1pt](A11)(A21)
\Edge[Direct,lw=1pt](A11)(A22)
\Edge[Direct,lw=1pt](A11)(A23)
\Edge[Direct,lw=1pt](A11)(A21)
\Edge[Direct,lw=1pt](A12)(A21)
\Edge[Direct,lw=1pt](A12)(A22)
\Edge[Direct,lw=1pt](A12)(A23)
\Edge[Direct,lw=1pt](A13)(A21)
\Edge[Direct,lw=1pt](A13)(A22)
\Edge[Direct,lw=1pt](A13)(A23)

\Edge[Direct,lw=1pt](A23)(Ah1)
\Edge[Direct,lw=1pt](A23)(Ah2)
\Edge[Direct,lw=1pt](A23)(Ah3)
\Edge[Direct,lw=1pt](A21)(Ah1)
\Edge[Direct,lw=1pt](A21)(Ah2)
\Edge[Direct,lw=1pt](A21)(Ah3)
\Edge[Direct,lw=1pt](A23)(Ah1)
\Edge[Direct,lw=1pt](A23)(Ah2)
\Edge[Direct,lw=1pt](A23)(Ah3)

\Edge[Direct,lw=1pt](Ah1)(A31)
\Edge[Direct,lw=1pt](Ah1)(A32)
\Edge[Direct,lw=1pt](Ah1)(A33)
\Edge[Direct,lw=1pt](Ah2)(A31)
\Edge[Direct,lw=1pt](Ah2)(A32)
\Edge[Direct,lw=1pt](Ah2)(A33)
\Edge[Direct,lw=1pt](Ah3)(A31)
\Edge[Direct,lw=1pt](Ah3)(A32)
\Edge[Direct,lw=1pt](Ah3)(A33)

\Edge[Direct,lw=1pt](A31)(tildeu)
\Edge[Direct,lw=1pt](A32)(tildeu)
\Edge[Direct,lw=1pt](A33)(tildeu)

\end{tikzpicture}
\caption{The architecture of a fully-connected feed-forward neural network. The trainable parameters are contained within the grey box.}
\label{figFCNN}
\end{center}
\end{figure}

A classical fully-connected, feed-forward neural network is described as a composition {\blue of} elementary functions. The NN contains $M$ layers, described as functions $L_i:\mathbb{R}^{N_i}\to\mathbb{R}^{N_{i+1}}$, where $N_{i+1}$ is the number of nodes in the layer. The layer is given by 
$$
L_i(x)=\varsigma_i(A_ix+b_i),
$$
where $\varsigma_i:\mathbb{R}\to\mathbb{R}$ is a user-prescribed {\it activation function}, acting componentwise, $A_i$ is an $N_{i+1}\times N_{i}$ matrix and $b_i\in\mathbb{R}^{N_{i+1}}$. The components of $A_i$ {\blue and} $b_i$ form the trainable parameters of the network. The fully-connected feed-forward NN is then given as the composition of the layers, as 
$$
u_{NN}(x)=(L_M\circ L_{M-1}\circ....\circ L_1)(x).
$$
The activation function of the last layer, $\varsigma_{M}$, is typically taken as the identity. Two common choices for the activation function in the remaining layers are the ReLU function, $\text{ReLU}(x)=\max(0,x)$ with $\text{ReLU}\in C(\mathbb{R})$ and differentiable everywhere but $x=0$, and the smooth hyperbolic tangent function $tanh$. {\blue\Cref{figFCNN} represents} the architecture schematically, where each circle corresponds to a node. The NN inherits its regularity from the activation function, so that the use of smooth activation functions in each layer {\blue implies} a smooth NN. 

\subsection{Approximation with classical NNs}
{\blue This section illustrates} via a 1D example the difficulty of Neural Networks (NNs) to approximate functions whose gradient is discontinuous. Let us consider the function $u\in H^1_0(0,\pi)$ given by 
\begin{equation}\label{1dsingularsol}
u(x)=\left\{ \begin{array}{ c c}
\sin(2x) & x\in \left[0,\frac{\pi}{2}\right],\\
\frac{1}{3}\sin(2x) & x\in \left(\frac{\pi}{2},\pi\right).
\end{array}\right.
\end{equation}
{\blue Note that} $u'$ is discontinuous at $x=\frac{\pi}{2}$. {\blue This target function corresponds to the weak solution of a 1D transmission problem that we will also consider in \Cref{subsecTrans1D}.}

We now approximate $u$ via a fully-connected feed-forward NN denoted $u_{NN}$. Explicitly, we consider three hidden layers of 20 neurons each and we compare the approximation with both $tanh$ and $\text{ReLU}$ activation functions. Whilst our focus is on using NNs to solve PDEs, we will use this example {\blue to} highlight the issues that may be present when training a NN by employing the simpler $H^1$ distance as a loss function. To do so,we approximate the $H^1$-norm of the error via Monte Carlo integration, as
\begin{equation}\label{Loss:1DMC}
\mathcal{L}(u_{NN})=\left(\frac{1}{N}\sum\limits_{i=1}^N|u(x_i)-u_{NN}(x_i)|^2+|u'(x_i)-u_{NN}'(x_i)|^2\right)^\frac{1}{2},
\end{equation}
where $\{x_i\}_{i=1}^N$ is a random uniform sample from $(0,\pi)$, and in each iteration we use a distinct sample. For our training, we use the Adam optimiser with an initial learning rate of $10^{-3}$, and train for 5,000 iterations. We take $N=2,500$ integration points.

\begin{figure}[H]
\begin{center}
\begin{subfigure}[b]{0.45\textwidth}
\includegraphics[height=0.55\textwidth]{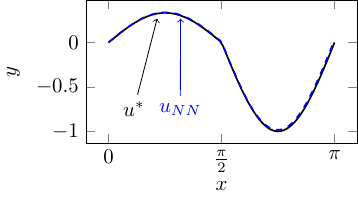}
\caption{Approximate and exact solution}
\end{subfigure}
\begin{subfigure}[b]{0.45\textwidth}
\includegraphics[height=0.55\textwidth]{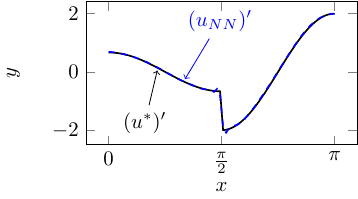}\caption{Approximate and exact gradient}
\end{subfigure}
\begin{subfigure}[b]{0.45\textwidth}
\includegraphics[height=0.55\textwidth]{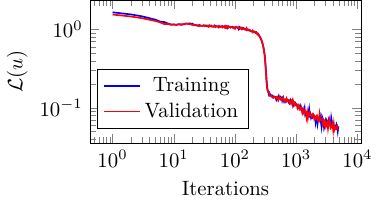}\caption{Loss evolution}
\end{subfigure}
\begin{subfigure}[b]{0.45\textwidth}
\includegraphics[height=0.55\textwidth]{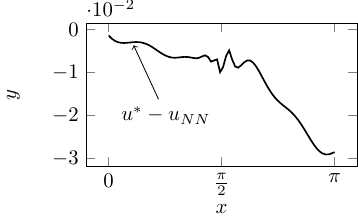}\caption{Approximation error}
\end{subfigure}
\begin{subfigure}[b]{0.45\textwidth}
\includegraphics[height=0.55\textwidth]{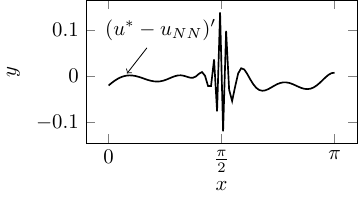}
\caption{Gradient approximation error}
\end{subfigure}
\caption{Classical architecture with Monte Carlo integration and tanh activation function.}
\label{Classical:tanh}
\end{center}
\end{figure}

Figure \ref{Classical:tanh} shows the approximation and its derivative when we take $tanh$ as the activation function. It also displays the errors in both the solution and the gradient, and the evolution of the loss functional. We observe that due to the Gibbs phenomenon, a uniform approximation of the discontinuous derivative is not possible as we are approximating $u$ employing smooth functions. Although a good approximation in $H^1$ should be possible, this would introduce high-frequency and order-one amplitude oscillations in the gradient near the discontinuity. We see in this case that the fluctuations are relatively localised, but still introduce significant errors in the $H^1$-norm. We obtain a relative error in $L^2$ of 2.6\% and relative $L^2$-error in the derivative of 4.9\%.

We expect these oscillations near the discontinuity to provoke numerical instabilities in the optimisation procedure. As we are using Monte Carlo integration, even if the oscillations are highly localised, it is likely that eventually a point will fall into a region where the error is high. Thus, producing large gradients in the loss with respect to the trainable parameters, even if the NN is a reasonable approximation of the target function in $H^1$.

\begin{figure}[H]
\begin{center}
\begin{subfigure}[b]{0.45\textwidth}
\includegraphics[height=0.55\textwidth]{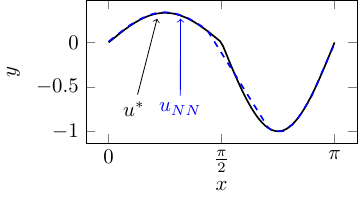}
\caption{Approximate and exact solution}
\end{subfigure}
\begin{subfigure}[b]{0.45\textwidth}
\includegraphics[height=0.55\textwidth]{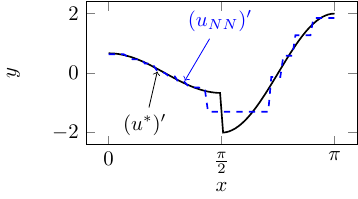}\caption{Approximate and exact gradient}
\end{subfigure}
\begin{subfigure}[b]{0.45\textwidth}
\includegraphics[height=0.55\textwidth]{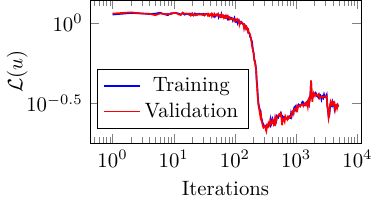}\caption{Loss evolution}
\end{subfigure}
\begin{subfigure}[b]{0.45\textwidth}
\includegraphics[height=0.55\textwidth]{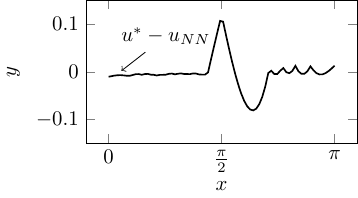}\caption{Approximation error}
\end{subfigure}
\begin{subfigure}[b]{0.45\textwidth}
\includegraphics[height=0.55\textwidth]{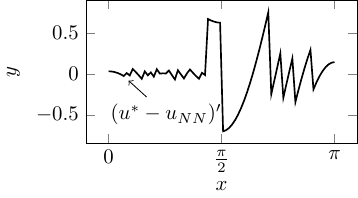}
\caption{Gradient approximation error}
\end{subfigure}
\caption{Classical architecture with Monte Carlo integration and ReLu activation function.}
\label{Classical:ReLU}
\end{center}
\end{figure}

Figure \ref{Classical:ReLU} shows the results of the same experiment when selecting $\text{ReLU}$  as the activation function. A fully-connected feed-forward NN with $\text{ReLU}$ activation is piecewise linear, admitting the same $H^1\setminus C^1$ regularity exhibited by the target function. However, in this case we observe an incredibly poor approximation, with visible instabilities in the loss. We obtain a relative error in $L^2$ of 6.4\% and a relative $L^2$ error in the derivative of 28.1\%. {\blue We argue that this appears because the loss contains spatial derivatives of the NN, so the gradient of the loss with respect to the trainable weights employs second derivatives of the activation function. As the second derivative of ReLU is a $\delta$-type function, the correct gradient cannot be captured by numerical quadrature, leading to poor optimisation properties.}

\subsection{ReCoNN for the 1D problem}

 Given the optimisation problems faced by classical NNs, we propose the incorporation of {\it a priori} knowledge of the regularity of the target function into the architecture. Let us assume that we know the nature and the location of the singularity in our solution. That is, we know that $u$ is $C^1$ except at $x=\frac{\pi}{2}$, $C^0$ on $(0,\pi)$, and admits a jump-singularity in $u'$ at $x=\frac{\pi}{2}$, whilst the height of the jump in the derivative is unknown. In this case, we consider a regularity-conforming architecture given by 
\begin{equation}\label{eqToyArchitecture}
u_{NN}(x)=\tilde{u}_0(x)+\tilde{u}_1(x)\frac{|x-\frac{\pi}{2}|}{2},
\end{equation}
where $\tilde{u}:\mathbb{R}\to\mathbb{R}^2$, $\tilde{u}=(\tilde{u}_0,\tilde{u}_1)$ and is a fully-connected feed-forward NN with smooth activation function. It is then immediate that $u_{NN}$ has the same regularity as the analytical solution and $u_{NN}'$ has left- and right-derivatives at $x=\frac{\pi}{2}$ given by 
\begin{equation}
\begin{split}\label{eqToyJump}
\lim\limits_{x\to\frac{\pi}{2}^\pm} u_{NN}'\left(x\right)= & \tilde{u}_0'\left(\frac{\pi}{2}\right)\pm\frac{1}{2}\tilde{u}_1\left(\frac{\pi}{2}\right).\\
\end{split}
\end{equation}

As both $\tilde{u}_0,\tilde{u}_1$ are outputs from a single NN $\tilde{u}:\mathbb{R}\to\mathbb{R}^2$ and $u_{NN}:\mathbb{R}\to\mathbb{R}$, the computational cost of such a parametrisation is marginal. Explicitly, using a fully connected feed-forward NN with {\blue three} hidden layers and 20 neurons per hidden layer, the classical architectures {\blue employ} 901 variables, and the conforming architecture {\blue incorporates} 922. Furthermore, whilst we expect the NN to have a greater approximation capability, we also observe that it allows to explicitly describe the discontinuity in the gradient via \eqref{eqToyJump}.  

\begin{figure}[H]
\begin{center}
\begin{subfigure}[b]{0.45\textwidth}
\includegraphics[height=0.55\textwidth]{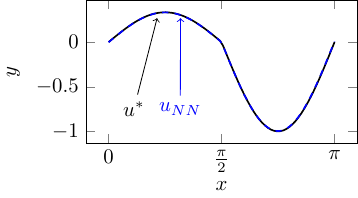}
\caption{Approximate and exact solution}
\end{subfigure}
\begin{subfigure}[b]{0.45\textwidth}
\includegraphics[height=0.55\textwidth]{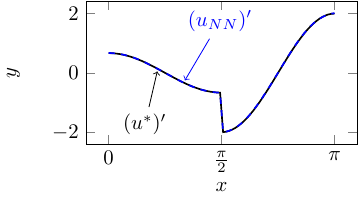}\caption{Approximate and exact gradient}
\end{subfigure}
\begin{subfigure}[b]{0.45\textwidth}
\includegraphics[height=0.55\textwidth]{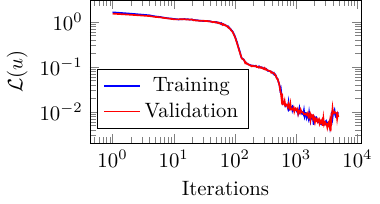}\caption{Loss evolution}
\end{subfigure}
\begin{subfigure}[b]{0.45\textwidth}
\includegraphics[height=0.55\textwidth]{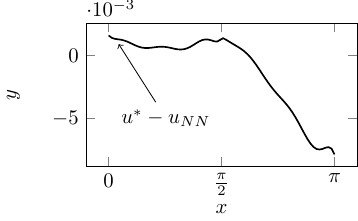}\caption{Approximation error}
\end{subfigure}
\begin{subfigure}[b]{0.45\textwidth}
\includegraphics[height=0.55\textwidth]{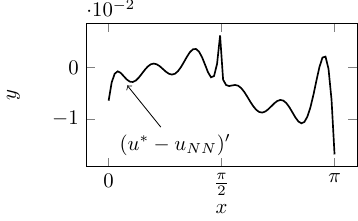}
\caption{Gradient approximation error}
\end{subfigure}
\caption{Regularity-conforming architecture with Monte Carlo integration and $tanh$ activation function.}
\label{ConformalArch}
\end{center}
\end{figure}

Figure \ref{ConformalArch} exhibits the results of minimizing the loss functional \eqref{Loss:1DMC} employing the regularity-conforming NN defined in \eqref{eqToyArchitecture} with $tanh$ activation function. We observe a far greater approximation capability: errors in both $u_{NN}$ and its derivative are uniformly small. We obtain a relative error of the solution in $L^2$ of 0.62\% and relative $L^2$ error in the derivative of 0.48\%. Furthermore, we observe far less oscillation in the loss during training, suggesting that the optimisation procedure is more stable. The instabilities at the end of training ocur when the loss is already small and the NN appears to have converged.

 From this simple example, we highlight the optimisation issues faced when approximating {\blue low-regularity functions} using classical NNs. We have observed that the incorporation of {\it a priori} knowledge of the regularity into the architecture greatly improves approximation capability and is partially explainable, in the sense that the singularity {\blue is described} according to \eqref{eqToyJump}. Our aim is to employ these ideas in the design of architectures for the 2D transmission problem. The key ingredient is the knowledge of the location and structure of the singularities, which we outline in the following section.

\section{Model Problem}\label{secModel}

We consider a bounded, polygonal domain $\Omega\subset\mathbb{R}^2$, and material data described by $\sigma:\Omega\to\mathbb{R}$. We assume that $\Omega$ consists of a finite number of non-intersecting subdomains, $\Omega_i$, and that $\sigma\in L^\infty(\Omega)$ is piecewise constant on the subdomains, so that $\sigma|_{\Omega_i}=\sigma_i$ for some $\sigma_i\in(0,\infty)$ and all $i=1,...,n$. {\blue We denote by $\Gamma\subset\Omega$ the set of discontinuities of $\sigma$ in $\Omega$.} Given $f\in L^2(\Omega)$, we consider the 2D transmission problem described by piecewise-constant material data and homogeneous Dirichlet boundary conditions given in weak form as: Find $u\in H^1_0(\Omega)$ such that 
\begin{equation}\label{eq:transmission}
\int_\Omega\sigma(x)\nabla u(x)\cdot\nabla v(x)+f(x)v(x)\,dx=0,
\end{equation}
for all $v\in H^1_0(\Omega)$. 

By considering test functions $v$ with $\text{supp}(v)\subset\Omega_i$, we see that $\Delta u=\frac{1}{\sigma_i}f$ weakly on each subdomain. {\blue Thus} by classical elliptic regularity \cite{bianca2023transmission} we have that $u\in H^2_{loc}(\Omega_i)$ for all $i$.
Furthermore, weak solutions satisfy $\sigma\nabla u\in H(\text{div};\Omega)$, and thus, for any interface $\Gamma\subset\Omega$ with normal vector $\nu$, we have that $\sigma\nabla u\cdot\nu$ is continuous across interfaces. Now, for a piecewise-smooth vector-valued function $v$, we define the jump operator at a point $x\in \Gamma$ by 
\begin{equation}
[v(x)]=\lim\limits_{t\to 0} \left(v(x+t\nu(x))-v(x-t\nu(x))\right)\cdot\nu(x).
\end{equation}

If we have sufficiently smooth solutions $u$, we obtain the following strong-form equation with interface conditions: 
\begin{equation}\label{StrongEqInter}
\begin{split}
\sigma_i\Delta u(x)=& f(x)\hspace{1cm}(x\in\Omega_i)\\
[\sigma(x)\nabla u(x)]=&0 \hspace{1cm}\left(x\in\bigcup\limits_{i=1}^n \partial\Omega_i\setminus\partial\Omega\right), \\
u(x)=&0 \hspace{1cm} (x\in\partial\Omega). 
\end{split}
\end{equation}

We aim to define architectures that can describe various singularities based on $\sigma$ and the geometry of $\Omega$.

{\blue If $\Omega$ is a convex set and the set of interfaces $\Gamma$ corresponds to non-intersecting, $C^1$, closed curves in $\Omega$, then the exact solution satisfies $u\in H^2(\Omega_i)$ for all $i$. In this case, the only singularities present are jump discontinuities in the normal component of the gradient, arising as a consequence of the continuity of the flux, $[\sigma(x)\nabla u(x)]=0$. 
}

\subsection{Gradient discontinuities with point singularities}\label{subsecInterface}
Besides the lack of regularity caused by jumps {\blue in the gradient} across interfaces, the transmission problem \eqref{eq:transmission} also admits power-like singularities at corners of the subdomains. In this work, {\blue we} restrict ourselves to the case of polygonal domains and consider two particular types of point singularities:
\begin{itemize}
\item Re-entrant corners, where the interior angle of the subdomain is above $180^o$. 
\item Material vertices, where three or more subdomains share a common vertex. 
\end{itemize}
 Other types {\blue of singularities} are possible, such as {\blue when} a material interface meets the boundary, or when the interface between two materials admits a corner. {\blue While both types of singularities may be considered with the methodology presented in this work, for simplicity we will avoid the details here.}

In each case, the solutions admit a decomposition into a regular part, which is in $H^2(\Omega_i)$, and a singular part \cite{blum1982finite}. Specifically, if $\{x_i\}_{i=1}^k$ index the corresponding vertices, we have that 
\begin{equation}\label{eq:decomp}
u(x)=w(x)+\sum\limits_{i=1}^k\sum\limits_{j=1}^{N_i}\kappa_{ij}s_{ij}(x),
\end{equation}
where $w\in H^2(\Omega_i)$ for all $i=1,...,n$, and $\kappa_{ij}$ are scalars known as \textit{stress intensity factors}. The singular functions $s_{ij}$ in \eqref{eq:decomp}, represented in polar coordinates centred at $x_i$, are of the form 
\begin{equation}\label{eq:singular}
s_{ij}(r,\theta)=\eta(r)r^{\lambda_{ij}}\phi_{ij}(\theta),
\end{equation}
where $\eta$ is an arbitrary, smooth cutoff function taking value $\eta(r)=1$ for $r<\delta_1$ and $\eta(r)=0$ for $r>\delta_2$, and $\phi_{ij}$ is an eigenvector of a Sturm-Liouville problem related to the geometry of the subdomains and has corresponding eigenvalue $\lambda^2_{ij}$. 

\subsubsection{Re-entrant corners}\label{subsubsecreen}
In the case of a polygonal domain that admits a re-entrant corner, {\blue such as the L-shape domain example of \Cref{figGeometryL},} for simplicity we take $\sigma=1$. It is known \cite{brenner1999multigrid} that the singular functions $s_i$ in \eqref{eq:singular} can be written in polar coordinates centred at $x_i$ as 
\begin{equation}
s_i(r,\theta)=\eta(r)r^{\lambda_i}\sin(\lambda_i(\theta-\omega_i)).
\end{equation}
$\lambda_i$ corresponds to the geometry so that $2\pi\lambda_i$ is the interior angle of the re-entrant corner, and $\omega_i$ is an appropriate phase-shift to ensure the Dirichlet condition.

\subsubsection{Vertices of material interfaces}\label{subsubsecmat}
\begin{figure}\begin{center}
\begin{subfigure}[b]{0.3\textwidth}
\begin{center}
\includegraphics[width=0.8\textwidth]{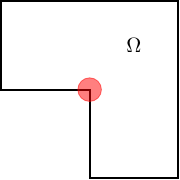}
\caption{Geometry of a re-entrant corner in the L-Shaped domain. The red circle highlights the re-entrant corner.}\label{figGeometryL}
\end{center}
\end{subfigure}\hspace{0.2\textwidth}
\begin{subfigure}[b]{0.3\textwidth}
\begin{center}
\includegraphics[width=0.8\textwidth]{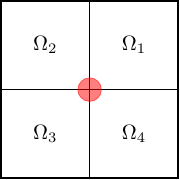}
\caption{Geometry of a domain with a vertex at material interfaces. The red circle highlights the singular point.}\label{figGeometry4}
\end{center}
\end{subfigure}
\caption{Examples of geometries that permit point singularities.}\end{center}
\end{figure}

In the case where three or more subdomains share a common vertex, we have that $w$ in the decomposition \eqref{eq:decomp} admits discontinuities in the gradient across interfaces. Also, the angular functions $\phi_{ij}$ in \eqref{eq:singular} are now periodic functions corresponding to eigenvectors of the weak Sturm-Liouville problem \cite{brenner2003multigrid}
\begin{equation}
\int_0^{2\pi}\sigma(\theta)\phi_{ij}'(\theta)v'(\theta)\,d\theta=\int_0^{2\pi}\sigma(\theta)\lambda_{ij}^2\phi_{ij}(\theta)v(\theta)\,d\theta
\end{equation}
for all periodic $v\in H^1(0,2\pi)$ and $\lambda_{ij}\in (0,1)$. Note that we have interpreted $\sigma$, via abuse of notation, as a function only of angle in a vecinity of the vertex $x_i$ {\blue as in  \Cref{figGeometry4}}.

Finally, as $\sigma$ admits discontinuities, this implies that $\phi_{ij}'$ can also admit discontinuities, although it will satisfy a continuous flux condition, i.e., $\sigma\phi_{ij}'$ is continuous.

\section{Architectures}\label{secArchitec}
In this section, we define the regularity-conforming architectures we consider for each of the cases presented in Section \ref{secModel}.

\subsection{Gradient discontinuities without point singularities}\label{subsecJump}
\begin{figure}\begin{center}
\includegraphics[width=0.3\textwidth]{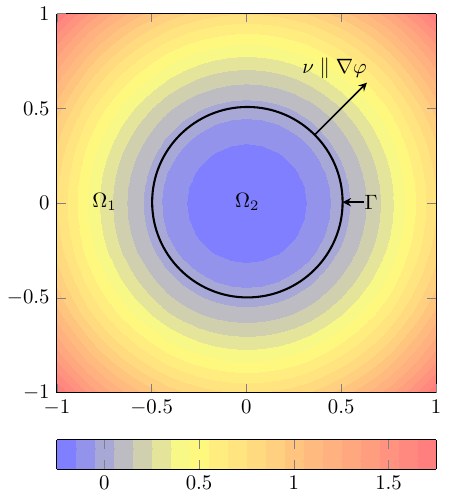}
\caption{Geometry of the domain $\Omega$ and the corresponding subdomains $\Omega_1,\Omega_2$ {\blue when $\Gamma$ is a smooth, closed curve}. The contourplot corresponds to the values of $\varphi(x)=|x|^2-0.25$.}\label{figGeometryInterface}\end{center}
\end{figure}

{\blue We turn to the definition of a ReCoNN in the case of smooth interfaces. We consider the simplified case where $\Omega$ is divided into two sub-domains $\Omega_1,\Omega_2$ separated by a smooth interface $\Gamma$, which forms a closed curve inside $\Omega$. As a modelling assumption, we assume that there exists a smooth function $\varphi:\Omega\to\mathbb{R}$ so that $\Omega_1=\{x\in\Omega:\varphi(x)>0\}$, $\Omega_2=\{x\in\Omega:\varphi(x)<0\}$, and $\Gamma$ corresponds to the zero set of $\varphi$. As a technical assumption, we assume that $|\nabla \varphi(x)|\neq 0$ for all $x\in\Gamma$. We remark that as $\Gamma$ is a level set of $\varphi$, we have that $\nu(x)$ is parallel to $\nabla \varphi(x)$ for all $x\in\Gamma$. In particular, we may write $\nabla \varphi(x)=|\nabla\varphi(x)|\nu(x)$. \Cref{figGeometryInterface} illustrates this geometry.

}

We consider the following architecture
\begin{equation}\label{eq:archinter}
u_{NN}(x)=\tilde{u}_0(x)+\tilde{u}_1(x)|\varphi(x)|,
\end{equation}
where $\tilde{u}(x)=(\tilde{u}_0(x),\tilde{u}_1(x))$, $\tilde{u}:\mathbb{R}^2\to\mathbb{R}^2$ is a fully-connected feed-forward NN with smooth activation function and $\varphi$ is a smooth function describing the interfaces introduced in Section \ref{subsecInterface}. This is a generalisation of \eqref{eqToyArchitecture} to interfaces in two dimensions.

As $\Gamma$ is a level set of $\varphi$, we have that $\nu(x)$, the normal vector to $\Gamma$ at a point $x\in\Gamma$, must be parallel to $\nabla\varphi(x)$. In particular, we may write $|\nabla\varphi(x)|\nu(x)=\nabla \varphi(x)$, by taking $\nu$ oriented so it points inwards to the domain where $\varphi$ is positive. In this case, we can calculate the gradient of $u_{NN}$ away from $\Gamma$ as 
\begin{equation}
\nabla u_{NN}(x)=\nabla \tilde{u}_0(x)+\nabla\tilde{u}_1(x)|\varphi(x)|+\text{sign}(\varphi(x))\tilde{u}_1(x)\nabla\varphi(x), 
\end{equation}
and for $x\in \Gamma$, recalling that $|\nabla\varphi(x)|\nu(x)=\nabla\varphi(x)$, we have that 
\begin{equation}\label{eq:grads}
\begin{split}
\lim\limits_{t\to 0^\pm}\nabla u_{NN}(x+t\nu(x)) =& \nabla \tilde{u}_0(x)\pm\tilde{u}_1(x)\nabla\varphi(x) = \nabla \tilde{u}_0(x)\pm\tilde{u}_1(x)|\nabla\varphi(x)|\nu(x).
\end{split}
\end{equation}
Thus, we may evaluate the left- and right-normal derivatives of $u_{NN}$ at the interface as
 \begin{equation}
\begin{split}
\lim\limits_{t\to 0^\pm}\nu(x)\cdot \nabla u_{NN}(x+t\nu(x)) = & \frac{\partial\tilde u_0}{\partial\nu} \pm \tilde{u}_1(x)|\nabla \varphi(x)|,
\end{split}
\end{equation}
and the jump in $[\nabla u_{NN}(x)]$ across $\Gamma$ is then calculated as $2|\nabla \varphi(x)|\tilde{u}_1(x)$. For this reason, as a technical condition, we require that $|\nabla\varphi(x)|\neq 0$ on $\Gamma$, else $u_{NN}$ will have continuous normal derivative across the interface.

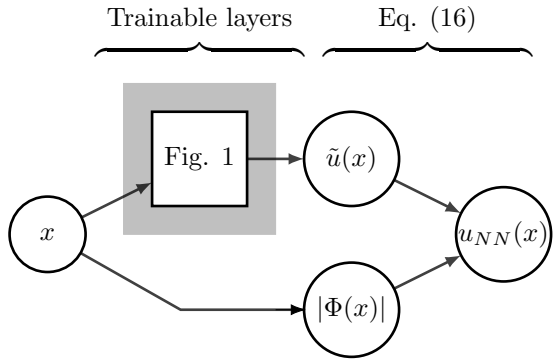
\begin{figure}[H]
\begin{center}
\begin{tikzpicture}

\filldraw [fill=black, draw=black,opacity=0.25] (1,-1) rectangle (3,1);

\Vertex[x=0,y=-1,label=$x$,color=white,size=1,fontsize=\normalsize]{X}

\Vertex[x=2,y=0,color=white,size=1.25,shape=rectangle,label=Fig. \ref{figFCNN},,fontsize=\normalsize]{NN1}
\Vertex[x=4,y=-2,label=$|\Phi(x)|$,color=white,size=1.25,fontsize=\normalsize]{vfx}
\Vertex[x=4,y=0,label=$\tilde{u}(x)$,color=white,size=1.25,fontsize=\normalsize]{tu1}
\Vertex[x=6,y=-1,label=$u_{NN}(x)$,color=white,size=1.25,fontsize=\normalsize]{outp}

\draw (2,1.85) node {Trainable layers};
\draw (2,1.5) node {$\overbrace{\hspace{2.75cm}}^\text{\hspace{1cm}}$};

\draw (5,1.85) node {Eq. \eqref{eq:archinterp}};
\draw (5,1.5) node {$\overbrace{\hspace{2.75cm}}^\text{\hspace{1cm}}$};

\Edge[Direct,lw=1pt](X)(NN1)
\Edge[Direct,lw=1pt](NN1)(tu1)
\Edge[Direct,lw=1pt,path = {X,{1.75,-2},vfx}](X)(vfx)
\Edge[Direct,lw=1pt](tu1)(outp)
\Edge[Direct,lw=1pt](vfx)(outp)

\draw[-{Latex[length=2mm]}] (3,-2) -- (3.4,-2);

\end{tikzpicture}

\caption{Representation of the ReCoNN used to approximate interfacial singularities. The node labelled ``Fig. \ref{figFCNN}" corresponds to a fully-connected feed-forward neural network, as in the cited figure. The assembly of $\tilde{u}$ and $|\Phi(x)|$ to yield $u_{NN}(x)$ is as in \eqref{eq:archinterp}. The trainable parameters are entirely within the fully-connected feed-forward network, indicated by the grey box.}
\label{figInterfaceReconn}
\end{center}
\end{figure}

We can easily generalise the architecture defined in \eqref{eq:archinter} to the case when the interfaces can be described as a union of smooth curves $\{\Gamma_p\}_{p=1}^P$, with each described as a zero-set of a smooth function $\{\varphi_p\}_{p=1}^P$. First, we define a feed-forward fully connected NN with smooth activation function $\tilde{u}=(\tilde{u}_0,\tilde{u}_1,...,\tilde{u}_P)$, $\tilde{u}:\mathbb{R}^2\to\mathbb{R}^{P+1}$, and then the regularity-conforming architecture as 
\begin{equation}\label{eq:archinterp}
u_{NN}(x)=\tilde{u}_0(x)+\sum\limits_{p=1}^P\tilde{u}_p(x)|\varphi_p(x)|. 
\end{equation}
Provided the curves admit intersections of zero length, $\varphi_{p}\neq 0$ on $\Gamma_{p'}$ for $p\neq p'$, similarly to \eqref{eq:grads}, the jumps in the gradient of $u_{NN}$ across $\Gamma_p$ are given by 

\begin{equation}
\begin{split}
\lim\limits_{t\to 0^{\pm}}\nabla u_{NN}(x+t\nu(x)) =& \nabla \left(\tilde{u}_0(x)+\sum\limits_{p'\neq p}\tilde{u}_{p'}(x)|\varphi_{p'}(x)|\right)\pm \tilde{u}_p(x)|\nabla\varphi_p(x)|\nu(x).
\end{split}
\end{equation}
Here, all gradients are well-defined away from the interfaces, and may be computed numerically via {\it autodiff}. As before, the normal component of the jump in the gradient across each interface $\Gamma_p$ is easily evaluated as 
\begin{equation}
[\nabla u_{NN}(x)]_p=2|\nabla \varphi_p(x)|\tilde{u}_p(x). 
\end{equation} 
Finally, we illustrate the architecture \eqref{eq:archinterp} in Figure \ref{figInterfaceReconn}, where $\Phi=(\varphi_1,...,\varphi_P)$, and the absolute value is taken componentwise. The grey box represents a fully-connected feed-forward NN with smooth activation function, containing the trainable parameters of the NN.

\subsection{Gradient discontinuities with point singularities}\label{subsecSingular}

Inspired by the regularity statement for re-entrant corners and material vertices, we define appropriate architectures in the next subsections.

\subsubsection{Re-entrant corners}\label{subsubsecCorner}

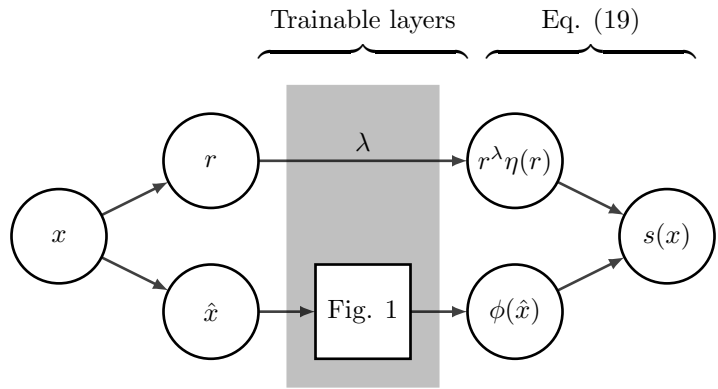
\begin{figure}[H]\begin{center}
\begin{tikzpicture}

\filldraw [fill=black, draw=black,opacity=0.25] (2.5,-0.5) rectangle (4.5,3.5);

\Vertex[x=-0.5,y=1.5,label=$x$,color=white,size=1.25,fontsize=\normalsize]{X}
\Vertex[x=1.5,y=2.5,color=white,size=1.25,fontsize=\normalsize,label=$r$]{xr}
\Vertex[x=1.5,y=0.5,label=$\hat{x}$,color=white,size=1.25,fontsize=\normalsize]{xhat}
\Vertex[x=3.5,y=0.5,color=white,size=1.25,shape=rectangle,label=Fig. \ref{figFCNN},fontsize=\normalsize]{NN1}
\Vertex[x=5.5,y=0.5,label=$\phi(\hat{x})$,color=white,size=1.25,fontsize=\normalsize]{pout}
\Vertex[x=7.5,y=1.5,label=$s(x)$,size=1.25,fontsize=\normalsize,color=white]{output}
\Vertex[x=5.5,y=2.5,label=$r^\lambda\eta(r)$,color=white,size=1.25,fontsize=\normalsize]{pr}

\Edge[Direct,lw=1pt](X)(xr)
\Edge[Direct,lw=1pt](X)(xhat)
\Edge[Direct,lw=1pt](xhat)(NN1)
\Edge[Direct,lw=1pt](NN1)(pout)
\Edge[Direct,lw=1pt](xr)(pr)
\Edge[Direct,lw=1pt](pr)(output)
\Edge[Direct,lw=1pt](pout)(output)

\draw (3.5,2.75) node {$\lambda$};
\draw (6.5,4.35) node {Eq. \eqref{eq:cornersingular}};
\draw (6.5,4) node {$\overbrace{\hspace{2.75cm}}^\text{\hspace{1cm}}$};

\draw (3.5,4.35) node {Trainable layers};
\draw (3.5,4) node {$\overbrace{\hspace{2.75cm}}^\text{\hspace{1cm}}$};

\end{tikzpicture}
\caption{Architecture of the singular component of the model. The rectangle labelled ``Fig. \ref{figFCNN}" corresponds to a fully connected feed forward network, as in the cited figure. The trainable parameters consist of $\lambda$ and the parameters of the fully-connected feed-forward sub-network. We have defined $r=|x-x_i|$ and $\hat{x}=\frac{1}{r}(x-x_i)$ for notational brevity.}\label{figSingularUnit}
\end{center}
\end{figure}

For the architecture corresponding to re-entrant corners presented in Section \ref{subsubsecreen}, we first define the singular component of the network. At each re-entrant corner, indexed by $i$, we consider  
\begin{equation}\label{eq:cornersingular}
s_i(x)=|x-x_i|^{\lambda_{i}}\eta(|x-x_i|)\phi_{i}\left(\frac{x-x_i}{|x-x_i|}\right).
\end{equation}
Here, $x_i$ corresponds to the vertex of the re-entrant corner, and $\phi_i:\mathbb{R}^2\to\mathbb{R}$ is a fully-connected feed-forward neural network. To ensure that $\phi_{i}$ depends only on the angle between $x$ and $x_i$, the unit vector in the direction {\blue of} $x-x_i$ {\blue is fed into the network}. Finally, $\eta$ is a user-prescribed {\blue and} sufficiently differentiable non-trainable cutoff function, and $\lambda_i$ is a trainable parameter. We illustrate the architecture in Figure \ref{figSingularUnit}.

Although it might be tempting to simply consider the architecture given by $u(x)=w(x)+\sum\limits_{i=1}^k s_i(x)$, with $w$ a classical feed-forward fully-connected NN, preliminary experiments revealed that the presence of the cutoff function $\eta$ introduced large gradients of $u$ into the interior of the domain. This led to a poorer approximation and slow convergence, as $w$ would have to compensate against these effects. To remedy this, we consider the architecture 
\begin{equation}\label{eq:archreen}
u_{NN}(x)=w_0(x)+\sum\limits_{i=1}^k \eta(|x-x_i|)w_i(x)+\sum\limits_{i=1}^k s_i(x),
\end{equation}
where $w=(w_0,...,w_{k}):\mathbb{R}^2\to\mathbb{R}^{k+1}$ is a fully-connected feed-forward neural network. Introducing the terms $\{w_i\}_{i=1}^k$ to compensate for the gradients in $\eta$ {\blue provides} better approximability and a more stable convergence.  

\begin{figure}[H]\begin{center}
\begin{tikzpicture}

\filldraw [fill=black, draw=black,opacity=0.25] (1.5,-1) rectangle (3.5,2.5);

\Vertex[x=0,y=0,color=white,size=1.25,label=$x$,fontsize=\normalsize]{X}
\Vertex[x=2.5,y=1.5,color=white,size=1.25,shape=rectangle,label=Fig. \ref{figFCNN} ,fontsize=\normalsize]{NN1}
\Vertex[x=2.5,y=0,color=white,size=1.25,shape =rectangle,label= Fig. \ref{figSingularUnit} ,fontsize=\normalsize]{NN2}
\Vertex[x=5,y=1.5,color=white,size=1.25,label=$w(x)$,fontsize=\normalsize]{w}
\Vertex[x=5,y=-1.5,color=white,size=1.25,label=$\eta(x)$,fontsize=\normalsize]{eta}
\Vertex[x=5,y=0,color=white,size=1.25,label=$s_1(x)$,fontsize=\normalsize]{s}
\Vertex[x=7,y=0,color=white,size=1.25,label=$u_{NN}(x)$,fontsize=\normalsize]{out}

\Edge[Direct,lw=1pt](X)(NN1)
\Edge[Direct,lw=1pt](X)(NN2)
\Edge[Direct,lw=1pt](NN1)(w)
\Edge[Direct,lw=1pt](NN2)(s)
\Edge[Direct,lw=1pt](w)(out)
\Edge[Direct,lw=1pt](s)(out)
\Edge[Direct,lw=1pt](eta)(out)
\Edge[Direct,lw=1pt,path = {X,{1.75,-1.5},eta}](X)(eta)
\draw[-{Latex[length=2mm]}] (2,-1.5) -- (4.4,-1.5);

\draw (2.5,3.15) node {Trainable Layers};
\draw (2.5,2.75) node {$\overbrace{\hspace{2.25cm}}^\text{\hspace{1cm}}$};

\draw (6,3.15) node {Eq. \eqref{eq:archreen}};
\draw (6,2.75) node {$\overbrace{\hspace{2.75cm}}^\text{\hspace{1cm}}$};
\end{tikzpicture}
\caption{Representation of the ReCoNN used to approximate re-entrant corners. All trainable parameters are within the sub-networks indicated by the gray box, with the sub-networks corresponding to the cited figures. }
\label{figArchReen}
\end{center}
\end{figure}
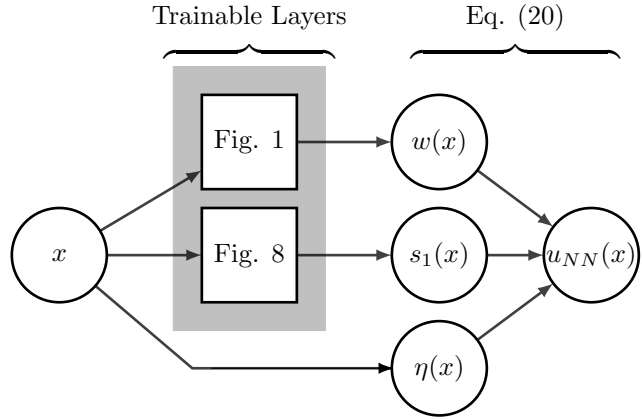

Figure \ref{figArchReen} illustrates architecture \ref{eq:archreen} in the simplified case of $k=1$. The trainable parameters are then the weights and biases of the fully-connected feed-forward NN that takes $\hat{x}\mapsto\phi(\hat{x})$, those of the NN that maps $x\mapsto w(x)$, and the scalar parameter $\lambda$.

This architecture is essentially two distinct neural networks in parallel, with the final network, $u_{NN}$, being {\blue their sum}. Unlike the case of Section \ref{subsecJump}, as the singular unit $s(x)$ and the $H^2$-component of $u_{NN}$ are based on distinct coordinate systems, it is not straightforward to implement a single network with the required properties. Nonetheless, the singular function effectively has a one-dimensional input and $\phi$ needs to approximate a ``simple” function (a sine in this case). Therefore, it is expected that a much smaller NN may be used for $\phi$, so it will not add a significant computational cost. The same observation applies to the following subsection.

\subsubsection{Vertices of material interfaces}\label{subsubsecVertices}

As discussed in Section \ref{subsubsecmat}, in the case of interior material vertices, the derivatives of the angular functions $\phi_{ij}$ can admit discontinuities. This may be rectified, however, by considering the singular functions to have an architecture analogous to that in \ref{subsecJump} and introducing discontinuities in its derivative to the architecture. 

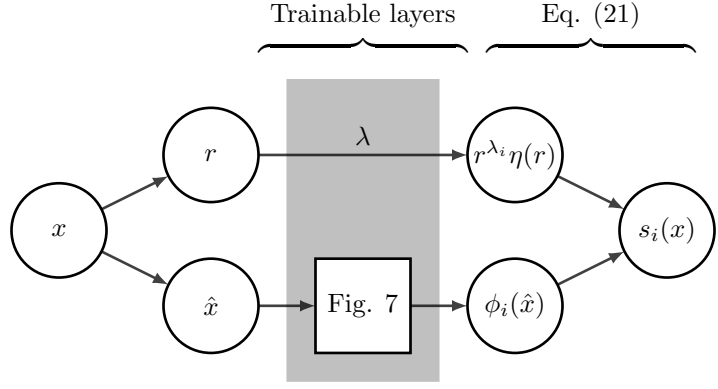
\begin{figure}[H]\begin{center}
\begin{tikzpicture}

\filldraw [fill=black, draw=black,opacity=0.25] (2.5,-0.5) rectangle (4.5,3.5);

\Vertex[x=-0.5,y=1.5,label=$x$,color=white,size=1.25,fontsize=\normalsize]{X}
\Vertex[x=1.5,y=2.5,color=white,size=1.25,fontsize=\normalsize,label=$r$]{xr}
\Vertex[x=1.5,y=0.5,label=$\hat{x}$,color=white,size=1.25,fontsize=\normalsize]{xhat}
\Vertex[x=3.5,y=0.5,color=white,size=1.25,shape=rectangle,label=Fig. \ref{figInterfaceReconn},fontsize=\normalsize]{NN1}
\Vertex[x=5.5,y=0.5,label=$\phi_i(\hat{x})$,color=white,size=1.25,fontsize=\normalsize]{pout}
\Vertex[x=7.5,y=1.5,label=$s_i(x)$,size=1.25,fontsize=\normalsize,color=white]{output}
\Vertex[x=5.5,y=2.5,label=$r^{\lambda_i}\eta(r)$,color=white,size=1.25,fontsize=\normalsize]{pr}

\Edge[Direct,lw=1pt](X)(xr)
\Edge[Direct,lw=1pt](X)(xhat)
\Edge[Direct,lw=1pt](xhat)(NN1)
\Edge[Direct,lw=1pt](NN1)(pout)
\Edge[Direct,lw=1pt](xr)(pr)
\Edge[Direct,lw=1pt](pr)(output)
\Edge[Direct,lw=1pt](pout)(output)

\draw (3.5,2.75) node {$\lambda$};
\draw (6.5,4.35) node {Eq. \eqref{MaterialSingular}};
\draw (6.5,4) node {$\overbrace{\hspace{2.75cm}}^\text{\hspace{1cm}}$};

\draw (3.5,4.35) node {Trainable layers};
\draw (3.5,4) node {$\overbrace{\hspace{2.75cm}}^\text{\hspace{1cm}}$};

\end{tikzpicture}
\caption{Architecture of the singular part of the model with material vertices. All trainable parameters are contained in the mappings indicated by the gray box. $r=|x-x_i|$ and $\hat{x}=\frac{1}{r}|x-x_i|$. }
\label{figSingInterface}
\end{center}
\end{figure}

First, we consider the architecture for the singular unit with interfaces. Heuristically, this is nothing more than a combination of the architectures proposed in Sections \ref{subsecJump} and \ref{subsubsecCorner}. As there may be several singular functions present at each material vertex, the user must consider a number $M_i$ of singular functions to include at each material vertex. Then, at each vertex, indexed by $i$ and centred at $x_i$, we define the singular components as 
\begin{equation}\label{MaterialSingular}
s_i(x)=\sum\limits_{j=1}^{M_i}\eta(|x-x_i|)|x-x_i|^{\lambda_{ij}}\phi_{ij}\left(\frac{x-x_i}{|x-x_i|}\right), 
\end{equation}
where $\phi_i=(\phi_{ij})_{j=1}^{M_i}:\mathbb{R}^2\to\mathbb{R}^{M_i}$ is a neural network whose architecture is given by \eqref{eq:archinterp}. The functions that define the interfaces as level sets and {\blue are} used to define the architecture in \eqref{eq:archinterp} may be described in polar coordinates centred at $x_i$. The parameters $\lambda_i=(\lambda_{ij})_{j=1}^{M_i}$ are trainable, whilst $\eta$ {\blue is again} a user-prescribed, sufficiently smooth, cutoff function. Figure \ref{figSingInterface} shows a schematic representation of the architecture in the simplified case when $M_i=1$.

Then, {\blue after defining the singular unit with interfaces,} we can include it in an architecture analogous to that of Figure \ref{figArchReen}. For simplicity, we {\blue consider} the case with a single singular vertex, with the case of multiple vertices being a straightforward extension. As previously, {\blue we build} $P$ smooth functions, $(\varphi_p)_{p=1}^P=\Phi$ {\blue whose zero sets are the given interfaces}. We then consider a fully-connected feed-forward neural network $w:\mathbb{R}^2\to\mathbb{R}^{2P+2}$ with components denoted $\tilde{w}(x)=((w_{i,p})_{p=0}^P)_{i=1}^2$. We then define
\begin{equation}\label{eq:wInteriorVertices}
w(x)=w_{1,0}(x)+w_{2,0}(x)\eta(|x-x_0|)+\sum\limits_{p=1}^P \Big(w_{1,p}(x)+w_{2,p}(x)\eta(|x-x_0|)\Big)|\varphi_p(x)|.
\end{equation}

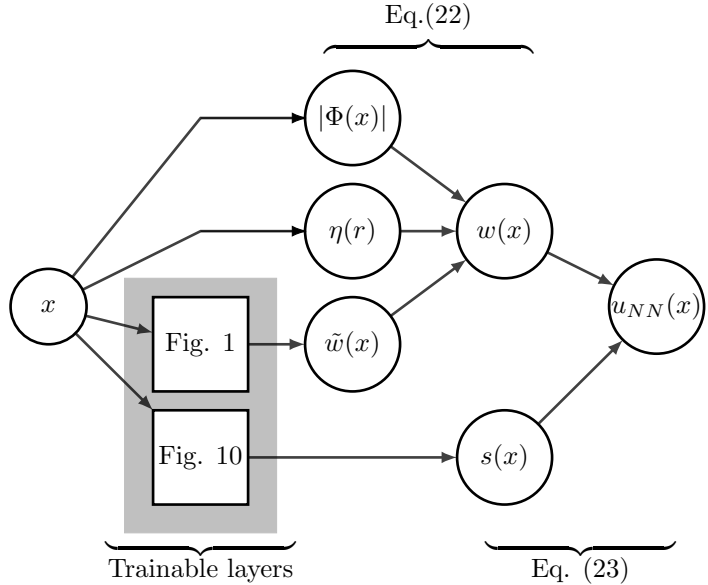
\begin{figure}[H]
\begin{center}
\begin{tikzpicture}

\filldraw [fill=black, draw=black,opacity=0.25] (1,-3) rectangle (3,0.375);

\Vertex[x=0,y=0,label=$x$,color=white,size=1,fontsize=\normalsize]{X}

\Vertex[x=2,y=-0.5,color=white,size=1.25,shape=rectangle, label=Fig. \ref{figFCNN},fontsize=\normalsize]{NN1}
\Vertex[x=4,y=-0.5,color=white,label=$\tilde{w}(x)$,size=1.25,fontsize=\normalsize]{w01}
\Vertex[x=4,y=2.5,label=$|\Phi(x)|$,color=white,size=1.25,fontsize=\normalsize]{vfx}
\Vertex[x=6,y=1,label=$w(x)$,color=white,size=1.25,fontsize=\normalsize]{w}
\Vertex[x=4,y=1,color=white,label=$\eta(r)$,size=1.25,fontsize=\normalsize]{eta}

\Vertex[x=2,y=-2,color=white,size=1.25,shape =rectangle,label=Fig. \ref{figSingInterface},fontsize=\normalsize]{NN2}
\Vertex[x=6,y=-2,color=white,size=1.25,label=$s(x)$,fontsize=\normalsize]{s}

\Vertex[x=8,y=0,color=white,size=1.25,label=$u_{NN}(x)$,fontsize=\normalsize]{out}

\Edge[Direct,lw=1pt](X)(NN1)
\Edge[Direct,lw=1pt](NN1)(w01)
\Edge[Direct,lw=1pt,path = {X,{2,2.5},vfx}](X)(vfx)
\Edge[Direct,lw=1pt](vfx)(w)
\Edge[Direct,lw=1pt](w01)(w)
\Edge[Direct,lw=1pt](X)(NN2)
\Edge[Direct,lw=1pt](NN2)(s)
\Edge[Direct,lw=1pt](s)(out)
\Edge[Direct,lw=1pt](w)(out)
\Edge[Direct,lw=1pt,path = {X,{2,1},eta}](X)(eta)
\Edge[Direct,lw=1pt](eta)(w)

\draw (5,3.85) node {Eq.\eqref{eq:wInteriorVertices}};
\draw (5,3.5) node {$\overbrace{\hspace{2.75cm}}^\text{\hspace{1cm}}$};

\draw (2,-3.15) node {$\underbrace{\hspace{2.5cm}}^\text{\hspace{1cm}}$};
\draw (2,-3.5) node {Trainable layers};

\draw (7,-3.15) node {$\underbrace{\hspace{2.5cm}}^\text{\hspace{1cm}}$};
\draw (7,-3.5) node {Eq. \eqref{eq:InteriorVerticesU}};

\draw[-{Latex[length=2mm]}] (2,2.5) -- (3.4,2.5);
\draw[-{Latex[length=2mm]}] (2,1) -- (3.4,1);
\end{tikzpicture}
\caption{Architecture for problems with polygonal piecewise materials and interior material vertices in the simplified case of $k=1$. The trainable parameters are contained within the sub-networks highlighted by the grey box and $r=|x-x_1|$.}
\label{figIntMat}
\end{center}
\end{figure}

The inclusion of the functions $w_{2,p}$ is to compensate for the large gradients away from the singularity in the singular unit, as in the case of Figure \ref{figArchReen}. Finally, $u_{NN}$ is defined as 
\begin{equation}\label{eq:InteriorVerticesU}
u_{NN}(x)=w(x)+\sum\limits_{i=1}^ks_i(x)
\end{equation} 
{\blue \Cref{figIntMat} illustrates the architecture,} taking $k=1$ for simplicity. 

\section{Numerical results}\label{secNumerics}
We now test the architectures proposed in Section \ref{secArchitec} in different transmission problems presenting different types of singularities. 

\subsection{1D Transmission problem}\label{subsecTrans1D}

We consider the following one-dimensional transmission problem: Find $u\in H^1_0(0,\pi)$
\begin{equation}
\int_0^{\pi}\sigma(x) u'(x)v'(x)-4\sin(2x)v(x)\,dx =0
\end{equation}
for all $v\in H^1_0(0,\pi)$, where 
\begin{equation}\label{1dtransmition}
\sigma(x)=\left\{ \begin{array}{ c c}
3 & x\in \left[0,\frac{\pi}{2}\right],\\
1 & x\in \left(\frac{\pi}{2},\pi\right).
\end{array}\right.
\end{equation}

The unique solution of the weak equation (\ref{1dtransmition}) given in \eqref{1dsingularsol}. {\blue The solution satisfies} the interface equation, 
\begin{equation}
\begin{split}
\sigma(x)u''(x)= & -4\sin(2x) \hspace{1cm}\left(x\in(0,\pi)\setminus\left\{\frac{\pi}{2}\right\}\right),\\
\lim\limits_{x\to \frac{\pi}{2}^+}\sigma(x)u'(x)=&\lim\limits_{x\to \frac{\pi}{2}^-}\sigma(x)u'(x).
\end{split}
\end{equation}

For this problem, we employ the architecture given in \eqref{eqToyArchitecture}, which is of the form described in Section \ref{subsecJump}. As the jump in the gradient is explicitly computable, we consider the following PINNs loss functional given by 
\begin{equation}\label{Loss1dtrans}
\begin{split}
\mathcal{L}_{PDE}(u_{NN})=&\frac{1}{N}\sum\limits_{i=1}^N|\sigma(x_i)u_{NN}''(x_i)+4\sin(2x_i)|^2,\\
\mathcal{L}_{int}(u_{NN})=&\left|\lim\limits_{x\to\frac{\pi}{2}^\pm} \sigma(x)u_{NN}'\left(x\right)-\lim\limits_{x\to\frac{\pi}{2}^-} \sigma(x)u_{NN}'\left(x\right)\right|^2,\\
\mathcal{L}_{bc}(u_{NN})=&\left(u_{NN}(0)^2+u_{NN}(\pi)^2\right),\\
\mathcal{L}(u_{NN})=& \alpha_1\mathcal{L}_{PDE}(u_{NN})^\frac{1}{2}+\alpha_2\mathcal{L}_{int}(u_{NN})^\frac{1}{2}+\alpha_3\mathcal{L}_{bc}(u_{NN})^\frac{1}{2}.
\end{split}
\end{equation}
where $\{x_i\}_{i=1}^N$ are random points from a uniform distribution in $(0,\pi)$ and $w_i$ are positive weights which, for simplicity, we take $\alpha_1=\alpha_2=\alpha_3=1$. The three terms in \eqref{Loss1dtrans} correspond to the ODE away from the interface, the interface condition, and Dirichlet condition, respectively. All of them can be calculated explicitly using {\it autodiff} via the representation \eqref{eqToyJump}.

\begin{figure}[H]
\begin{center}
\begin{subfigure}[b]{0.45\textwidth}
\includegraphics[height=0.55\textwidth]{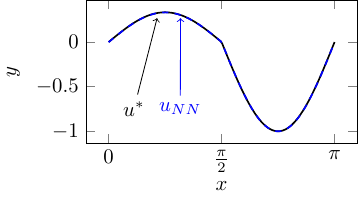}
\caption{Approximate and exact solution}
\end{subfigure}
\begin{subfigure}[b]{0.45\textwidth}
\includegraphics[height=0.55\textwidth]{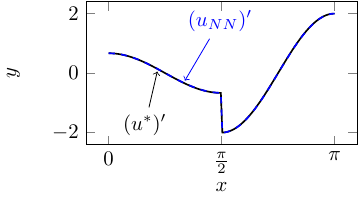}\caption{Approximate and exact gradient}
\end{subfigure}
\begin{subfigure}[b]{0.45\textwidth}
\includegraphics[height=0.55\textwidth]{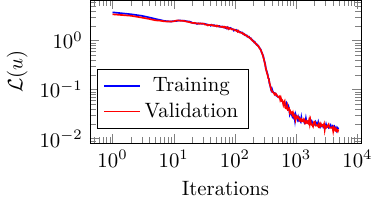}\caption{Loss evolution}
\end{subfigure}
\begin{subfigure}[b]{0.45\textwidth}
\includegraphics[height=0.55\textwidth]{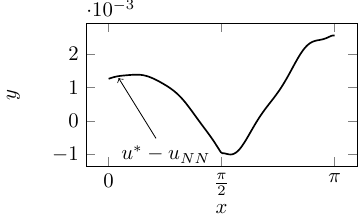}\caption{Approximation error}
\end{subfigure}
\begin{subfigure}[b]{0.45\textwidth}
\includegraphics[height=0.55\textwidth]{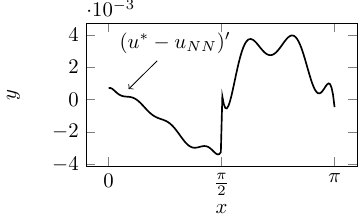}
\caption{Gradient approximation error}
\end{subfigure}
\caption{Regularity-conforming architecture with PINNs loss functional. We obtain a relative $L^2$-error of $0.25\%$ and a relative $L^2$-error in the derivative of $0.22\%$.}
\label{ConformalArchPINNs}
\end{center}
\end{figure}

We perform a similar optimisation strategy as in Section \ref{secMotivation}, taking 5000 iterations, $N=2500$, and an Adam optimiser with initial learning rate of $10^{-3}$. Figure \ref{ConformalArchPINNs} displays the approximation results and the evolution of the loss.  We obtain a relative error in $L^2$ of 0.25\% and {\blue a} relative $L^2$ error in the derivative of 0.22\%. We observe a good approximation in $H^1$ norm, and a less oscillatory behaviour of the loss during training comparing to the results in Section \ref{secMotivation}. As the jump condition is included in the loss, we see very good approximation of the gradient across the discontinuity.

\subsection{Smooth material interface}\label{subsecSmooth2D}

We consider now a 2D example {\blue with a} smooth material interface. Let us consider $\Omega=(-1,1)^2$, and $\sigma:\Omega\to\mathbb{R}$ to be given by 
\begin{equation}
\sigma(x)=\left\{\begin{array}{c c}
1 & |x|<\frac{1}{2},\\
3 & \text{else}. 
\end{array}\right.
\end{equation}
In this case, we {\blue express the interface} as the zero set of the smooth function $\varphi(x)=|x|^2-\frac{1}{4}$.
We consider the manufactured solution given by 
\begin{equation}
u^*(x)=\frac{(x_1-1)^2(x_2-1)^2(4x_1^2+4x_2^2-1)^2}{\sigma(x)},
\end{equation}
and see that for $f=\sigma\Delta u^*$, which is defined everywhere but {\blue on} the interface $\Gamma$, then $u^*$ satisfies the interface equation corresponding to $\sigma$, 
\begin{equation}\label{eqNoFlux}
\left[\sigma \nabla u^* \right]=0. 
\end{equation}
As the left- and right-normal derivatives can be numerically computed, the right-hand term of the continuous flux condition on the interface \eqref{eqNoFlux} can be readily computed for a candidate function $u_{NN}$. 

For the implementation we consider the regularity-conforming architecture \eqref{eq:archinter}, employing a fully-connected feed-forward neural network of {\blue three} hidden layers of 30 neurons each with $tanh$ activation function for $\tilde{u}:\mathbb{R}^2\to\mathbb{R}^2$. The PINNs loss function is the following:
\begin{equation}\label{LossInterf}
\begin{split}
\mathcal{L}_{PDE}(u_{NN})=&\frac{1}{N_1}\sum\limits_{i=1}^{N_1}|\sigma(x_i^1)\Delta u_{NN}(x_i^1)-f(x_i^1)|^2,\\
\mathcal{L}_{int}(u_{NN})=&\frac{1}{N_2}\sum\limits_{i=1}^{N_2}\left[\sigma(x^{2}_i) \nabla u_{NN}(x^2_i)\right]^2,\\
\mathcal{L}_{bc}(u_{NN})=&\frac{1}{N_3}\sum\limits_{i=1}^{N_3}|u_{NN}(x_i^3)|^2,\\
\mathcal{L}(u_{NN})=& \alpha_1\mathcal{L}_{PDE}(u_{NN})^\frac{1}{2}+\alpha_2\mathcal{L}_{int}(u_{NN})^\frac{1}{2}+\alpha_3\mathcal{L}_{bc}(u_{NN})^\frac{1}{2}.
\end{split}
\end{equation}
$\{x_i^1\}_{i=1}^{N_1},\{x_i^2\}_{i=1}^{N_2}$ and $\{x_i^3\}_{i=1}^{N_3}$ are random, uniformly sampled points at each iteration in $\Omega$, $\Gamma$, and $\partial\Omega$, respectively, and we select $N_1=N_2=N_3=1,000$. {\blue Employing the square root of each component} of the loss ensures better convergence when near a minimum, {\blue as noted in \cite{mishra2022estimates}}. We choose the weights $\alpha_1=1$, $\alpha_2=\sqrt{10}$, and $\alpha_3=10$, following \cite{wang2022and} where authors suggest that the weights for lower-order terms should be progressively higher to ensure good convergence towards a minimiser.

\begin{figure}[H]\begin{center}
\begin{subfigure}[h]{0.45\textwidth}\begin{center}
\includegraphics[height=0.85\textwidth]{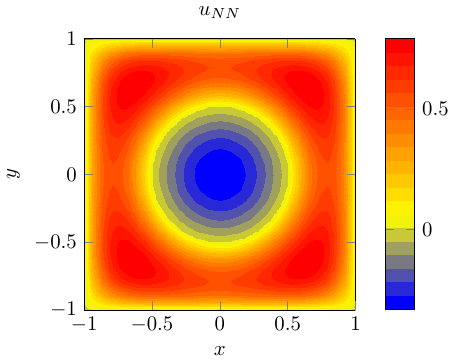}
\caption{Approximate solution $u_{NN}$}
\end{center}\end{subfigure}
\begin{subfigure}[h]{0.45\textwidth}\begin{center}
\includegraphics[height=0.85\textwidth]{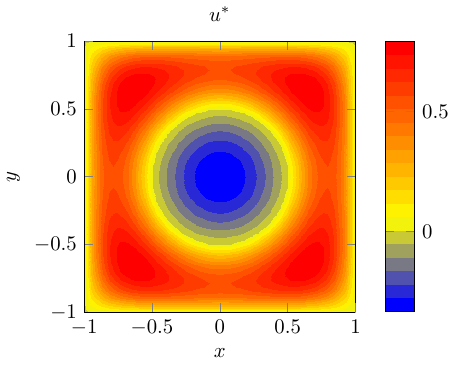}
\caption{Exact solution $u^*$}
\end{center}\end{subfigure}
\begin{subfigure}[h]{0.45\textwidth}\begin{center}
\includegraphics[height=0.85\textwidth]{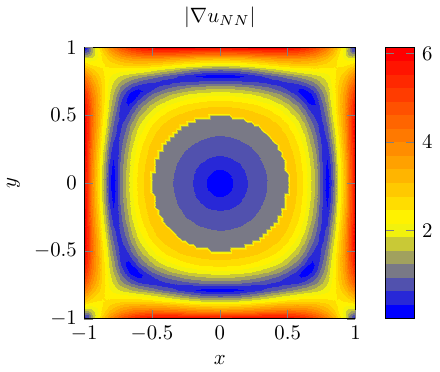}
\caption{Norm of the gradient of the approximate solution $|\nabla u_{NN}|$}
\end{center}\end{subfigure}
\begin{subfigure}[h]{0.45\textwidth}\begin{center}
\includegraphics[height=0.85\textwidth]{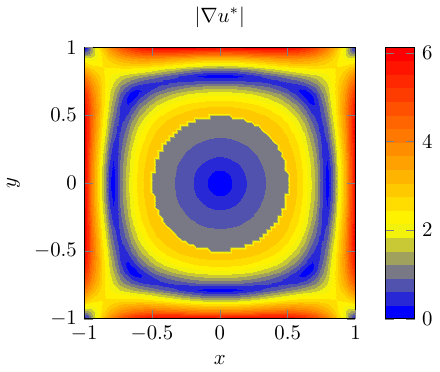}
\caption{Norm of the gradient of the exact solution $|\nabla u^*|$\\
$\left. \right.$}
\end{center}\end{subfigure}
\begin{subfigure}[h]{0.45\textwidth}\begin{center}
\includegraphics[height=0.85\textwidth]{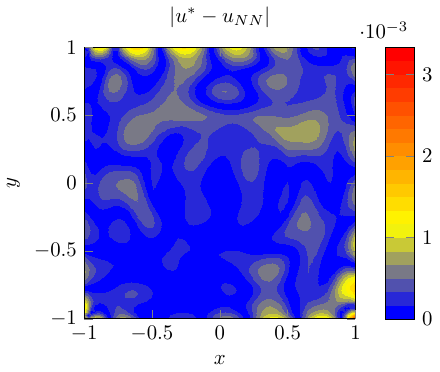}
\caption{Error of the solution $|u^*-u_{NN}|$}
\end{center}\end{subfigure}
\begin{subfigure}[h]{0.45\textwidth}\begin{center}
\includegraphics[height=0.85\textwidth]{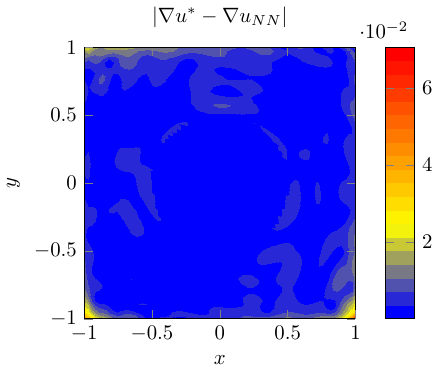}
\caption{Error of the gradient of the solution $|\nabla u^*-\nabla u_{NN}|$}
\end{center}\end{subfigure}
\caption{Exact and approximate solutions. We obtain an $L^2$-relative error of $0.16\%$ and an $L^2$-relative error in the gradient of $0.075\%$.}
\label{fig:SolutionInterf}\end{center}
\end{figure}

\begin{figure}[H]\begin{center}
\begin{subfigure}[b]{0.45\textwidth}\begin{center}
\includegraphics[width=0.95\textwidth]{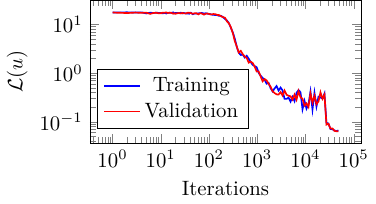}
\caption{Total loss}\end{center}\end{subfigure}
\begin{subfigure}[b]{0.45\textwidth}\begin{center}
\includegraphics[width=0.95\textwidth]{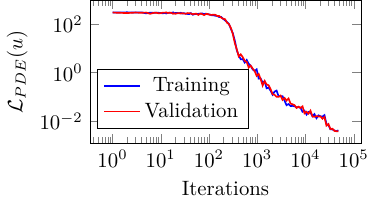}
\caption{PDE loss}\end{center}\end{subfigure}
\begin{subfigure}[b]{0.45\textwidth}\begin{center}
\includegraphics[width=0.95\textwidth]{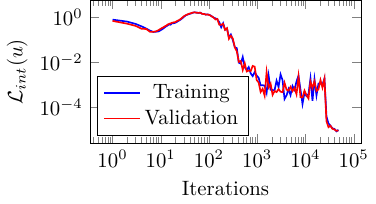}
\caption{Interface condition loss}\end{center}\end{subfigure}
\begin{subfigure}[b]{0.45\textwidth}\begin{center}
\includegraphics[width=0.95\textwidth]{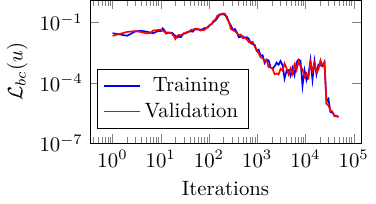}
\caption{Boundary condition loss}\end{center}\end{subfigure}
\caption{Loss evolution during training. }
\label{fig:LossInterf}
\end{center}
\end{figure}

For our optimisation strategy, we employ the Adam optimiser with 50,000 iterations. We set a learning rate of $10^{-3}$ for the first 25,000 iterations, and for the latter 25,000 iterations we reduce the learning rate via an exponential decay so that it is $10^{-6}$ at the end of training. Figure \ref{fig:SolutionInterf} displays the analytical and the approximated solutions, their gradients, and the errors. Figure \ref{fig:LossInterf} shows the evolution of the loss.

We observe good, uniform approximation of both the solution and its gradient, with relative $L^2$ errors of the solution and gradient of 0.16\% and 0.076\%, respectively. There is no visible Gibbs phenomenon present, with the errors in both $u$ and $\nabla u$ being more prominent towards the boundary of the domain, rather than the ring corresponding to the discontinuity in $\sigma$. The presence of errors on the boundary suggests that these particular errors may be reduced by a more appropriate weighting of the boundary term in the loss \eqref{LossInterf} or a more appropriate integration strategy in its evaluation. Nonetheless, our focus in this work is on the description of discontinuities, with these results showing strong approximation capabilities of our architecture with the PINNs-type loss including the interface condition.

\subsection{L-shape domain}\label{subsecLshape}
Let us now consider the L-shaped domain, $\Omega=(-1,1)^2\setminus[-1,0]^2$ with $\sigma(x)=1$ on the entirety of $\Omega$. This domain possesses a re-entrant corner as defined in Section \ref{subsubsecreen}, which provokes a singularity, and the corresponding singular function is given explicitly as
\begin{equation}
s_0(r,\theta)=r^{\frac{2}{3}}\sin\left(\frac{2}{3}\left(\theta-\frac{\pi}{2}\right)\right),
\end{equation}

We consider the source $f=\Delta u^*$ for the following constructed solution 
\begin{equation}
u^*(x)=s_0(x)(x_1^2-1)(x_2^2-1). 
\end{equation}
The regularity statement requires a sufficiently smooth cutoff function which may be chosen by the user. For our implementation, we first define the piecewise following polynomial 
\begin{equation}
\eta_0(t)=\left\{\begin{array}{ll}
 1& t<0,\\
 -6t^5 + 15t^4 - 10t^3 + 1 & 0\leq t\leq 1,\\
 0 & t>1,
\end{array}\right.
\end{equation}
which defines a globally $C^2$ function, equal to $1$ for $t<0$ and $0$ for $t>1$. We then define $\eta(r)=\eta_0\left(\frac{r}{\delta_2-\delta_1}+\frac{\delta_1}{\delta_1-\delta_2}\right)$, to give a globally $C^2$ cutoff function equal to $1$ for $r<\delta_1$ and $0$ for $r>\delta_2$. We set $\delta_2=0.5$ and $\delta_1=0.9$ in our implementation. 

We consider the PINN loss given by 
\begin{equation}\label{LossLshape}
\begin{split}
\mathcal{L}_{PDE}(u_{NN})=&\frac{1}{N_1}\sum\limits_{i=1}^{N_1}|\Delta u_{NN}(x_i^1)-f(x_i^1)|^2\omega(x_i^1),\\
\mathcal{L}_{bc}(u_{NN})=&\frac{1}{N_2}\sum\limits_{i=1}^{N_2}| u_{NN}(x_i^{2})|^2,\\
\mathcal{L}_{bc_\phi}(u_{NN})=& \phi(-1,0)^2+\phi(0,-1)^2,\\
\mathcal{L}(u_{NN})=&\alpha_1\mathcal{L}_{PDE}(u_{NN})^\frac{1}{2}+\alpha_2\mathcal{L}_{bc}(u_{NN})^\frac{1}{2}+\alpha_3\mathcal{L}_{bc_\phi}(u_{NN})^\frac{1}{2},\\
\end{split}
\end{equation}
The components of the loss \eqref{LossLshape} correspond to the PDE, the Dirichlet boundary condition on $u_{NN}$, and the Dirichlet boundary condition on the singular function $\phi$, respectively. Here, $\{x^1_i\}_{i=1}^{N_1}$ are randomly and uniformly sampled {\blue points} from $\Omega$ and $\{x^2_i\}_{i=1}^{N_2}$ from $\partial\Omega$. We select $N_1=N_2=1,000$ and, similarly to the interface problem, we chose higher weights for lower order terms by taking $\alpha_1=1$, $\alpha_2=10$. We take $\alpha_3=1$.

We highlight that the PDE component of the loss contains a weight function $\omega$, which we define as
\begin{equation}\label{Lshapeweight}
\omega(x)=\min(40|x|^2,1),
\end{equation}
so that $\mathcal{L}_{PDE}$ corresponds to a discretisation of 
\begin{equation}\label{LshapeweightPDE}
\int_\Omega |\Delta u_{NN}(x)-f(x)|^2\omega(x)\,dx.
\end{equation}
At an arbitrary candidate solution $u_{NN}$ that does not solve the PDE, an elementary calculation shows that the singular part of $u_{NN}$ has a Laplacian squared that behaves as $|\Delta u_{NN}|^2\sim r^{2\lambda-4}$. For $\lambda\in (0,1)$, this is not integrable in $\mathbb{R}^2$, and thus may introduce significant numerical instabilities. To avoid the issues that this would necessarily provoke, we employ the weight function $\omega$ \eqref{Lshapeweight}, which behaves as $ r^2$ near $x=0$ and is equal to $1$ for $|x|^2>0.025$. By introducing this weight, the integrand in \ref{LshapeweightPDE} scales, at worst, as $r^{2\lambda-2}$, and is thus integrable as $\lambda>0$. As $\omega$ is taken to be positive almost everywhere, the integral \ref{LshapeweightPDE} is zero if and only if $u_{NN}$ is a solution of the PDE almost everywhere. 

The use of this weighting function is inspired by regularity statements for transmission problems in weighted-$H^2$ spaces of a similar form introduced in \cite{li2010analysis}.  Here, the authors employ norms of the form $\int_\Omega |\nabla^2 u|^2\omega(x)\,dx$ for weight functions that scale as $r^\beta$ for $x$ near a singular point.

\begin{figure}[H]\begin{center}
\begin{subfigure}[b]{0.45\textwidth}\begin{center}
\includegraphics[height=0.85\textwidth]{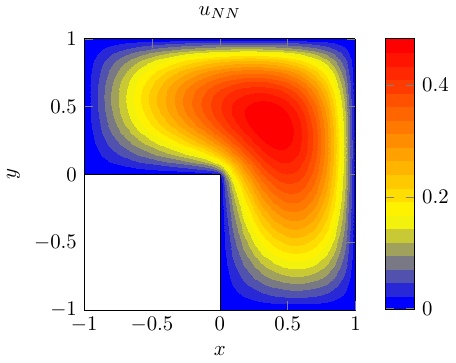}
\caption{Approximate solution $u_{NN}$.}
\end{center}\end{subfigure}
\begin{subfigure}[b]{0.45\textwidth}\begin{center}
\includegraphics[height=0.85\textwidth]{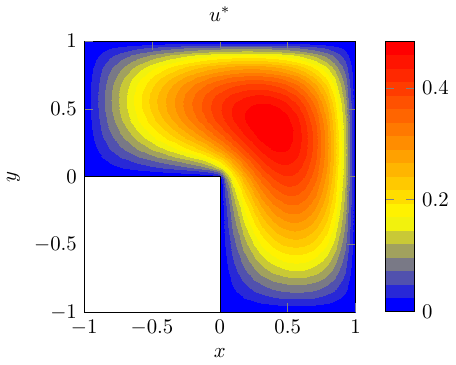}
\caption{Exact solution $u^*$.}
\end{center}\end{subfigure}
\begin{subfigure}[b]{0.45\textwidth}\begin{center}
\includegraphics[height=0.85\textwidth]{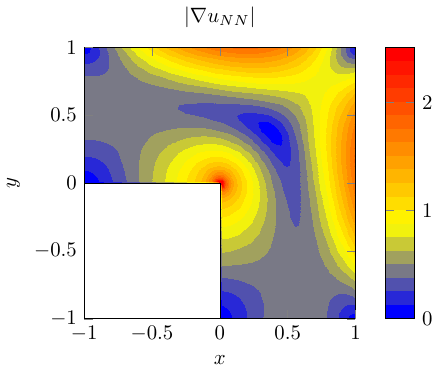}
\caption{Norm of the gradient of the approximate solution $|\nabla u_{NN}|$}
\end{center}\end{subfigure}
\begin{subfigure}[b]{0.45\textwidth}\begin{center}
\includegraphics[height=0.85\textwidth]{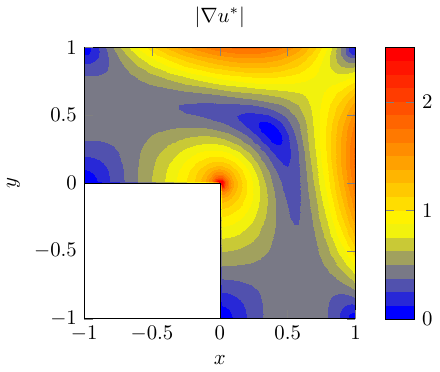}
\caption{Norm of the gradient of the exact solution $|\nabla u^*|$\\
$\left. \right.$}
\end{center}\end{subfigure}
\begin{subfigure}[b]{0.45\textwidth}\begin{center}
\includegraphics[height=0.85\textwidth]{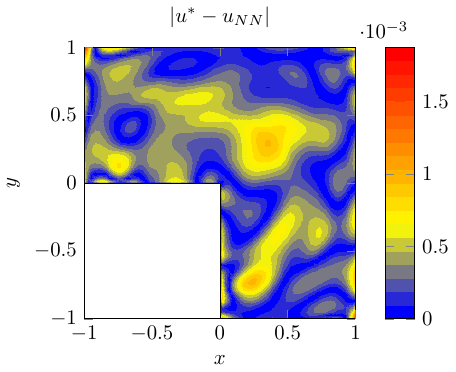}
\caption{Error of the solution $|u^*-u_{NN}|$.}
\end{center}\end{subfigure}
\begin{subfigure}[b]{0.45\textwidth}\begin{center}
\includegraphics[height=0.85\textwidth]{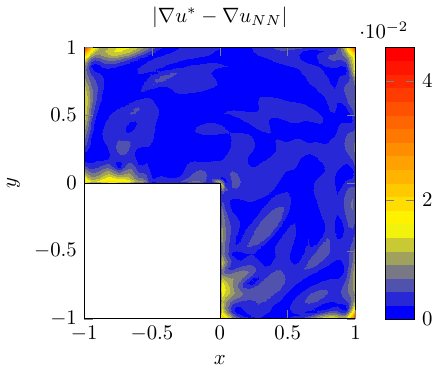}
\caption{Error of the gradient of the solution $|\nabla u^*-\nabla u_{NN}|$}
\end{center}\end{subfigure}
\caption{Exact and approximate solutions. We obtain a relative $L^2$-error of $0.15\%$ and a relative $L^2$-error in the gradient of $0.55\%$.}
\label{fig:Lshapesolutions}\end{center}
\end{figure}

\begin{figure}[H]\begin{center}
\begin{subfigure}[b]{0.45\textwidth}\begin{center}
\includegraphics[width=0.95\textwidth]{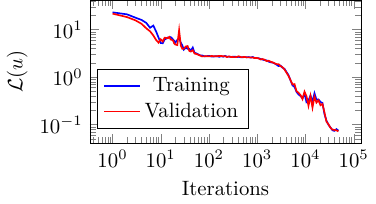}
\caption{Total loss}\end{center}\end{subfigure}
\begin{subfigure}[b]{0.45\textwidth}\begin{center}
\includegraphics[width=0.95\textwidth]{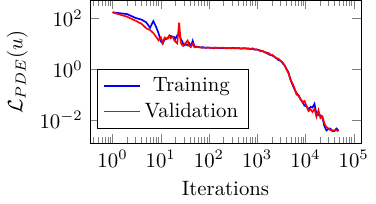}
\caption{PDE loss}\end{center}\end{subfigure}
\begin{subfigure}[b]{0.45\textwidth}\begin{center}
\includegraphics[width=0.95\textwidth]{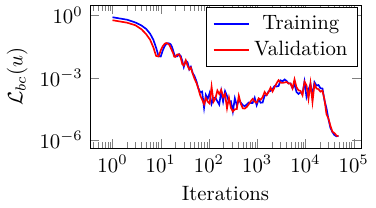}
\caption{Boundary condition loss}\end{center}\end{subfigure}
\begin{subfigure}[b]{0.45\textwidth}\begin{center}
\includegraphics[width=0.95\textwidth]{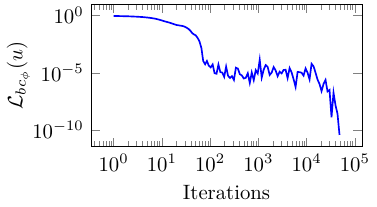}
\caption{Boundary condition loss}\end{center}\end{subfigure}
\caption{Loss evolution during training.}
\label{fig:LshapesoLoss}
\end{center}
\end{figure}

\begin{figure}[H]\begin{center}
\begin{subfigure}[t]{0.45\textwidth}\begin{center}
\includegraphics[height=0.6\textwidth]{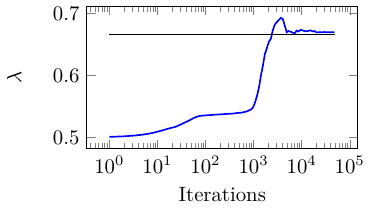}
\caption{Evolution of the exponent $\lambda$ during training and the exact value. }\end{center}\end{subfigure}
\begin{subfigure}[t]{0.45\textwidth}\begin{center}
\includegraphics[height=0.6\textwidth]{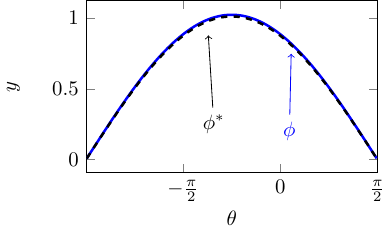}
\caption{The exact and approximate singular solutions as a function of angle.}\end{center}\end{subfigure}
\begin{subfigure}[t]{0.45\textwidth}\begin{center}
\includegraphics[height=0.6\textwidth]{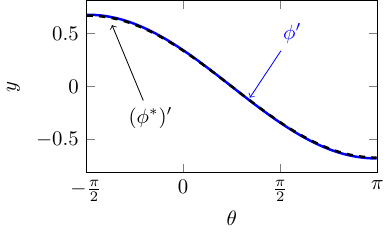}
\caption{The flux of the exact and approximate singular solutions as a function of angle.}\end{center}\end{subfigure}
\caption{Training and solutions of the singular function.}
\label{fig:Lshapesingular}\end{center}
\end{figure}

For the numerical results, we consider the regularity conforming architecture defined in \eqref{eq:archreen} with $k=1$. For that, we select $w:\mathbb{R}^2\to\mathbb{R}^2$ to be a fully-connected feed-forward NN with {\it tanh} activation function, three hidden layers and 30 neurons per hidden layer, resulting in 2012 trainable variables. We take $\phi$ to be a fully-connected feed-forward network with 3 hidden layers and 15 neurons per hidden layer, giving 542 trainable parameters. This gives a total of 2555 variables. We perform 50,000 iterations with the Adam optimiser, employing a learning rate of $10^{-3}$ for the first 25,000 iterations, then exponentially decreasing the learning rate to a final value of $10^{-6}$ for the last 25,000 iterations. Figure \ref{fig:Lshapesolutions} displays the approximation, its gradient, and the error; Figure \ref{fig:LshapesoLoss} shows the evolution of the loss; and Figure \ref{fig:Lshapesingular} the approximation of the singular function defined in \eqref{eq:cornersingular}.

We observe good approximation properties, both of the smooth component of the solution and the singular function, achieving a relative $L^2$-error of 0.15\% and a relative $L^2$-error in the gradient of 0.55\%. At the end of training, we obtain and exponent $\lambda=0.669$, which is a good approximation to the exact value $\lambda^*=\frac{2}{3}$. Whilst the error in the gradient at the origin is unbounded, as will always be the case unless the singular function and exponent are exact, the graphics show a good uniform approximation throughout the domain, with no significant errors present near the reentrant corner.  We see that both the singular solution and its derivative approximate well the exact solution, and the exponent $\lambda$ approximates the exact exponent sufficiently well as to avoid large, localised errors near the reentrant corner. Finally, as in Section \ref{subsecInterface}, we see larger errors in the gradient localised on the boundary, particularly in the corners.

\begin{figure}[H]\begin{center}
\begin{subfigure}[b]{0.45\textwidth}\begin{center}
\includegraphics[height=0.85\textwidth]{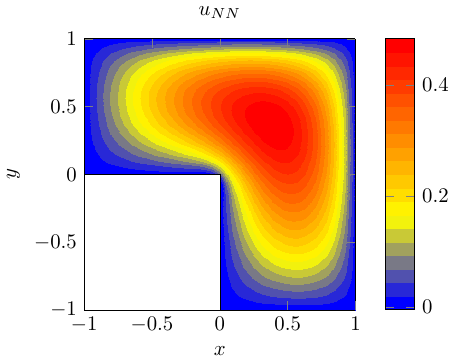}
\caption{Approximate solution $u_{NN}$.}
\end{center}\end{subfigure}
\begin{subfigure}[b]{0.45\textwidth}\begin{center}
\includegraphics[height=0.85\textwidth]{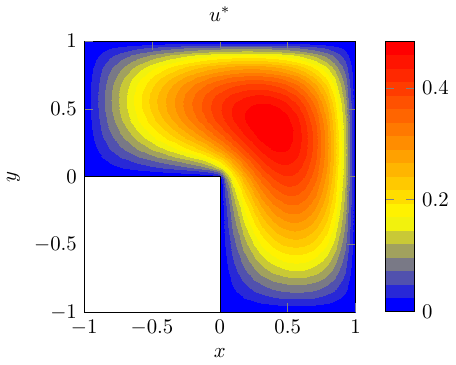}
\caption{Exact solution $u^*$.}
\end{center}\end{subfigure}
\begin{subfigure}[b]{0.45\textwidth}\begin{center}
\includegraphics[height=0.85\textwidth]{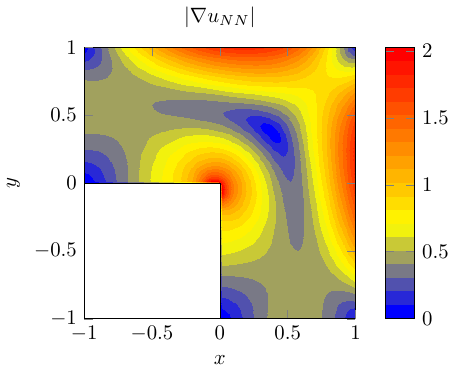}
\caption{Norm of the gradient of the approximate solution $|\nabla u_{NN}|$}
\end{center}\end{subfigure}
\begin{subfigure}[b]{0.45\textwidth}\begin{center}
\includegraphics[height=0.85\textwidth]{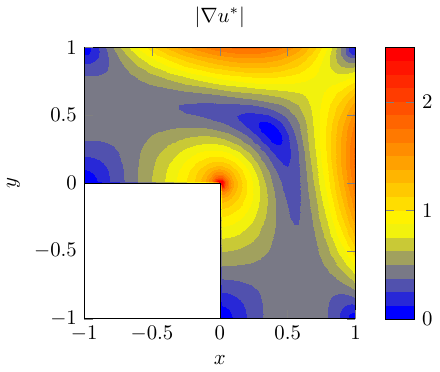}
\caption{Norm of the gradient of the exact solution $|\nabla u^*|$\\
$\left. \right.$ }
\end{center}\end{subfigure}
\begin{subfigure}[b]{0.45\textwidth}\begin{center}
\includegraphics[height=0.85\textwidth]{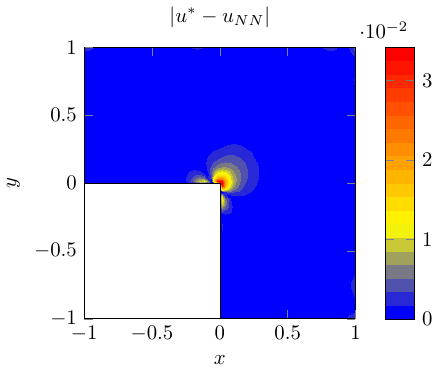}
\caption{Error of the solution $|u^*-u_{NN}|$.}
\end{center}\end{subfigure}
\begin{subfigure}[b]{0.45\textwidth}\begin{center}
\includegraphics[height=0.85\textwidth]{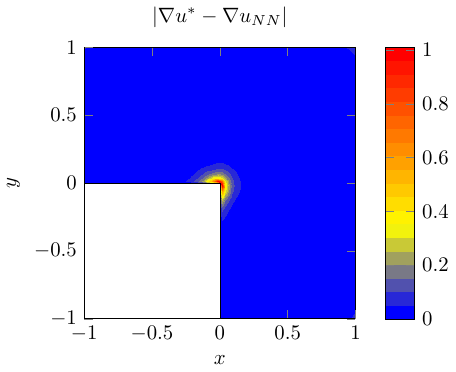}
\caption{Error of the gradient of the solution $|\nabla u^*-\nabla u_{NN}|$}
\end{center}\end{subfigure}
\caption{Exact and approximate solutions (Classical architecture). We obtain a relative $L^2$-error of $0.7\%$ and a relative $L^2$-error in the gradient of $5.7\%$.}
\label{fig:LshapesolutionsClass}\end{center}
\end{figure}

\begin{figure}[H]\begin{center}
\begin{subfigure}[b]{0.45\textwidth}\begin{center}
\includegraphics[width=0.95\textwidth]{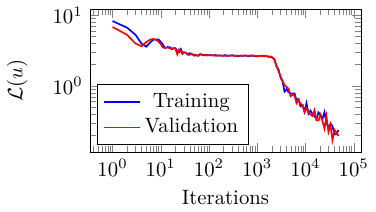}
\caption{Total loss}\end{center}\end{subfigure}
\begin{subfigure}[b]{0.45\textwidth}\begin{center}
\includegraphics[width=0.95\textwidth]{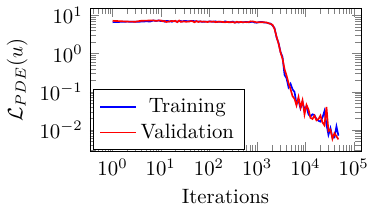}
\caption{PDE loss}\end{center}\end{subfigure}
\begin{subfigure}[b]{0.45\textwidth}\begin{center}
\includegraphics[width=0.95\textwidth]{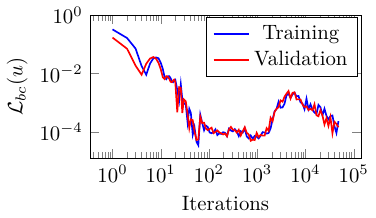}
\caption{Boundary condition loss}\end{center}\end{subfigure}
\caption{Loss evolution during training (Classical architecture).}
\label{fig:LshapeLossClass}
\end{center}
\end{figure}

For comparison, as $u^*$ satisfies the PDE in the classical sense on $\Omega$, that is, its Laplacian is well-defined and continuous on $\Omega$, we include a PINN implementation employing a classical architecture. We select a fully connected feed-forward NN with $tanh$ activation function, using three hidden layers and 35 nodes per hidden layer, resulting in marginally more trainable variables -- a total of 2661-- than the regularity-conforming architecture. We employ the same loss and training method as with the regularity-conforming network, {\blue but we exclude $\mathcal{L}_{bc_\phi}$} as there is no singular function in the classical architecture. Figures \ref{fig:LshapesolutionsClass} and  \ref{fig:LshapeLossClass} show the approximate solution and the {\blue loss evolution.}


In this case, we observe a relatively good approximation of $u$ itself, corresponding to a relative $L^2$-error of 0.7\%. Nonetheless, there is a highly localised error in the origin, leading to a relative $L^2$-error roughly four times greater than the regularity-conforming network. The failure of approximation of the gradient at the origin is even more pronounced, leading to order-one errors at the resolution of the plot, and producing relative $L^2$-errors in the gradient of 5.7\%. Curiously, we have that $\mathcal{L}_{PDE}$ at the end of training is comparable in the case of the classical and regularity-conforming architectures, whilst $\mathcal{L}_{bc}$ is two orders of magnitude higher. This indicates that the failure of the ability to approximate the singularity would appear to lead to non-compliance of the boundary condition rather than non-compliance of the PDE, resulting in errors that propagate into the interior. 

\subsection{Interior material vertices}\label{subsecVertices}
As a final example, we consider a problem that exhibits the singularities described in Section \ref{subsubsecmat}. We consider the domain $\Omega=(-1,1)^2$, where $\sigma$ is piecewise constant in each quadrant, i.e., 
\begin{equation}\label{SigmaMat}
\sigma(x)=\left\{\begin{array}{c c}
1 & x_1>0,x_2>0,\\
2 & x_1<0,x_2>0,\\
3 & x_1<0,x_2<0,\\
4 & x_1>0,x_2<0. 
\end{array}\right.
\end{equation}
 There exists only one exact singular solution, which can be written as $s^*(r,\theta)=r^\lambda \left(a_i\sin(\lambda\theta)+b_i\cos(\lambda\theta)\right)$, for $(r\cos(\theta),r\sin(\theta))\in\Omega_i$, $i=1,\ldots,4$. Numerically, we can find the coefficients to be $\lambda\approx 0.8599$, 
\begin{equation}
\begin{array}{l l l l}
a_1 = 3.584 & a_2 = 3.285 & a_3 = 2.474 & a_4 = 2.115\\
b_1 = -2.003 & b_2 = -0.6678 & b_3 = -1.0495 & b_4 = -0.5861.
\end{array}
\end{equation}

We consider the constructed solution given by 
$$u^*(x)=\cos\left(\frac{x_1\pi}{2}\right)\cos\left(\frac{x_2\pi}{2}\right)s^*(x),$$ 
and take $f=\sigma \nabla u^*$ so it satisfies the interface equation (\ref{StrongEqInter}) with $\sigma$ in \ref{SigmaMat}. The details of the construction of this solution are outlined in \Cref{App}. In this particular geometry, we consider the interfaces to be described by two functions: $\varphi_1(x)=x_1$ and $\varphi_2(x)=x_2$. We employ the regularity-conforming neural network architecture defined in Section \ref{subsubsecVertices} with the $C^2$ cutoff function defined in Section \ref{subsecLshape}.

As a loss function, we consider the following:
\begin{equation}
\begin{split}
\mathcal{L}_{PDE}(u_{NN})=&\frac{1}{N_1}\sum\limits_{i=1}^{N_1}|\sigma(x_i^1)\Delta u_{NN}(x_i^1)-f(x_i^1)|^2\omega_1(x^1_i),\\
\mathcal{L}_{int}(u_{NN})=&\frac{1}{N_2}\sum\limits_{i=1}^{N_2}\left[\sigma\nabla u_{NN}(x_i^{2})\right]^2\omega_2(x_i^2),\\
\mathcal{L}_{bc}(u_{NN})=&\frac{1}{N_3}\sum\limits_{i=1}^{N_3}| u_{NN}(x_i^{3})|^2,\\
\mathcal{L}_{bc_\phi}(u_{NN})=&\frac{1}{4}\sum\limits_{i=1}^4 \left[\sigma(x^4_i)\nabla\phi(x^4_i)\right]^2,\\
\mathcal{L}(u_{NN})=&\alpha_1\mathcal{L}_{PDE}(u_{NN})^\frac{1}{2}+\alpha_2\mathcal{L}_{int}(u_{NN})^\frac{1}{2}+\alpha_3\mathcal{L}_{bc}(u_{NN})^\frac{1}{2}+\alpha_4\mathcal{L}_{bc_\phi}(u_{NN})^\frac{1}{2}.\\
\end{split}
\end{equation}
Here, $\{x_i^1\}_{i=1}^{N_1}$,$\{x_i^2\}_{i=1}^{N_2}$, and $\{x_i^3\}_{i=1}^{N_3}$ are randomly, uniformly sampled from $\Omega$, $\Gamma$, and $\partial\Omega$, respectively, via a stratified sample that produces the same number of points in each quadrant. The term $\mathcal{L}_{bc_\phi}$ corresponds to the continuous flux condition on the singular function, and corresponds to summation over the four angles where the derivative of the singular function admits jumps. As in Section \ref{subsecLshape}, we employ the weight function $\omega_1$ to stabilise the integration of the singularity in the Laplacian, once again taking $\omega_1(x)=\min(40|x|^2,1)$. Similarly, the term $\mathcal{L}_{int}(u)$ scales as $r^{2\lambda-2}$, and corresponds to the discretisation of the  line integral $\displaystyle{\int_\Gamma [\sigma \nabla u]^2w_2(x)\,dx}$, therefore, it will not be integrable for a general trial function. Thus, we employ another weighting function $\omega_2(x)=\min(40|x|,1)$ to avoid numerical instabilities. Finally, we consider the weights to be $\alpha_1=1,\alpha_2=\sqrt{10},\alpha_3=10$ and $\alpha_4=1$, and we select $N_1=N_2=N_3=10^3$. 

\begin{figure}[H]\begin{center}
\begin{subfigure}[b]{0.45\textwidth}\begin{center}
\includegraphics[height=0.85\textwidth]{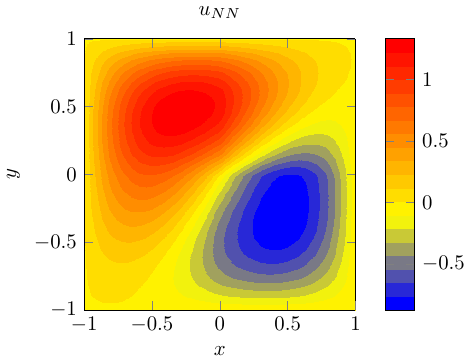}
\caption{Approximate solution $u_{NN}$.}
\end{center}\end{subfigure}
\begin{subfigure}[b]{0.45\textwidth}\begin{center}
\includegraphics[height=0.85\textwidth]{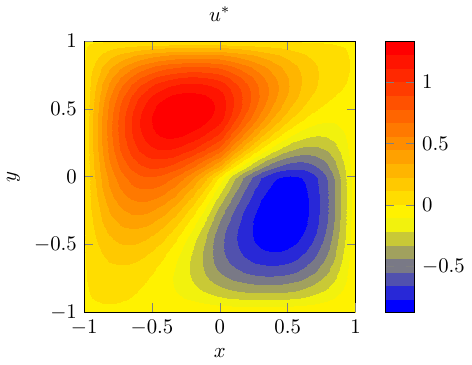}
\caption{Exact solution $u^*$.}
\end{center}\end{subfigure}
\begin{subfigure}[b]{0.45\textwidth}\begin{center}
\includegraphics[height=0.85\textwidth]{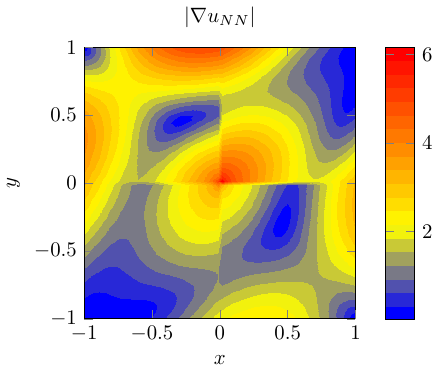}
\caption{Norm of the gradient of the approximate solution $|\nabla u_{NN}|$}
\end{center}\end{subfigure}
\begin{subfigure}[b]{0.45\textwidth}\begin{center}
\includegraphics[height=0.85\textwidth]{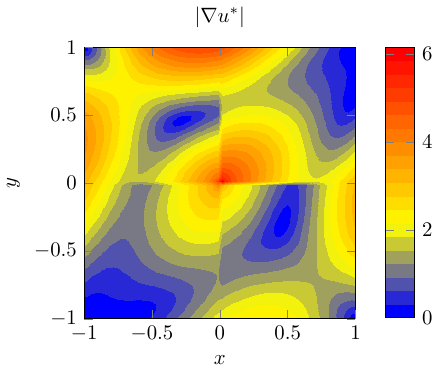}
\caption{Norm of the gradient of the exact solution $|\nabla u^*|$\\
$\left. \right.$}
\end{center}\end{subfigure}
\begin{subfigure}[b]{0.45\textwidth}\begin{center}
\includegraphics[height=0.85\textwidth]{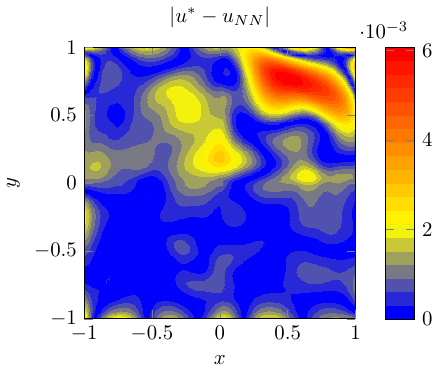}
\caption{Error of the solution $|u^*-u_{NN}|$.}
\end{center}\end{subfigure}
\begin{subfigure}[b]{0.45\textwidth}\begin{center}
\includegraphics[height=0.85\textwidth]{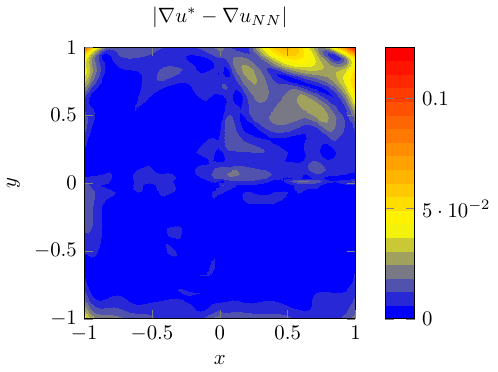}
\caption{Error of the gradient of the solution $|\nabla u^*-\nabla u_{NN}|$}
\end{center}\end{subfigure}
\caption{Exact and approximate solutions. We obtain a relative $L^2$-error of $0.25\%$ and a relative $L^2$-error in the gradient of $0.59\%$.}
\label{fig:MaterialSolutions}\end{center}
\end{figure}

\begin{figure}[H]\begin{center}
\begin{subfigure}[b]{0.45\textwidth}\begin{center}
\includegraphics[width=0.95\textwidth]{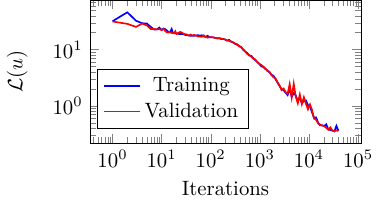}
\caption{Total loss}\end{center}\end{subfigure}
\begin{subfigure}[b]{0.45\textwidth}\begin{center}
\includegraphics[width=0.95\textwidth]{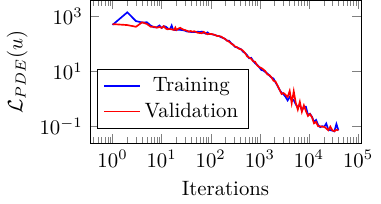}
\caption{PDE loss}\end{center}\end{subfigure}
\begin{subfigure}[b]{0.45\textwidth}\begin{center}
\includegraphics[width=0.95\textwidth]{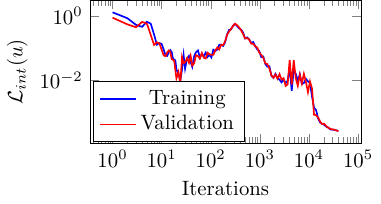}
\caption{Interface condition loss}\end{center}\end{subfigure}
\begin{subfigure}[b]{0.45\textwidth}\begin{center}
\includegraphics[width=0.95\textwidth]{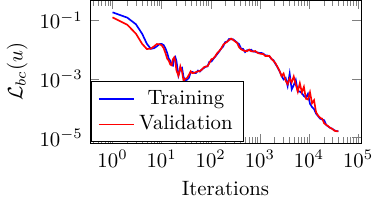}
\caption{Boundary condition loss}\end{center}\end{subfigure}
\caption{Loss evolution during training.}
\label{fig:MaterialLoss}
\end{center}
\end{figure}

\begin{figure}[H]\begin{center}
\begin{subfigure}[t]{0.45\textwidth}\begin{center}
\includegraphics[width=0.95\textwidth]{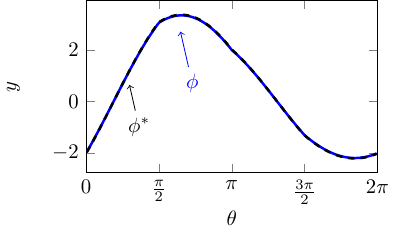}
\caption{The exact and approximate singular solutions as a function of angle.}\end{center}\end{subfigure}
\begin{subfigure}[t]{0.45\textwidth}\begin{center}
\includegraphics[width=0.95\textwidth]{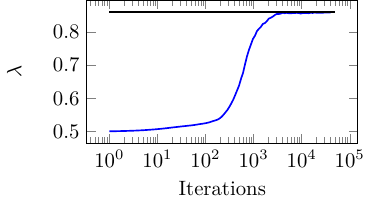}
\caption{Evolution of the exponent $\lambda$ during training and the exact value.}\end{center}\end{subfigure}
\begin{subfigure}[t]{0.45\textwidth}\begin{center}
\includegraphics[width=0.95\textwidth]{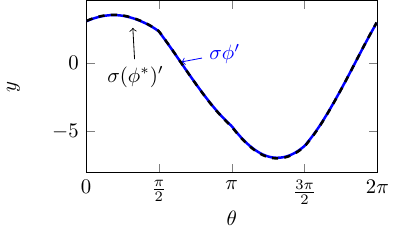}
\caption{The flux of the exact and approximate singular solutions as a function of angle.}\end{center}\end{subfigure}
\begin{subfigure}[t]{0.45\textwidth}\begin{center}
\includegraphics[width=0.95\textwidth]{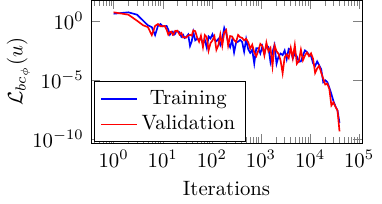}
\caption{Loss corresponding to the flux of the singular function during training.}\end{center}\end{subfigure}
\caption{Training and solutions of the singular function.}
\label{fig:MaterialSingular}\end{center}
\end{figure}

As before, we train with an Adam optimiser for 50,000 iterations with a learning rate of $10^{-3}$ for the first 25,000 iterations and an exponential decay to $10^{-6}$ for the last 25,000 iterations. For the architecture, we consider three hidden layers of 30 neurons each in the regular part of the architecture, and three hidden layers of 15 neurons each in the singular unit with interfaces. Along with the trainable parameter that describes the exponent, this gives a total of 2710 trainable variables. Figure \ref{fig:MaterialSolutions} displays the approximation, its gradient, and the error, and  Figure \ref{fig:MaterialLoss} shows the evolution of the loss. Figure \ref{fig:MaterialSingular} shows the approximation of the singular function defined in \eqref{MaterialSingular} and the term in the loss corresponding to it.

We observe good approximations of both the solution an the gradient, and uniformly, we obtain a $L^2$-relative error of 0.25\% and a $L^2$-relative error in the gradient of 0.59\%. We also see a good approximation of the singular function and its flux, as well as the approximation of the exponent as $\lambda=0.8591$, in comparison to its exact value of $\lambda^*=0.8599$. {\blue As in} the case of corner singularities, any  error in the singular function or exponent will lead to unbounded error in the gradient. Nonetheless, whilst limited by the resolution of the graphs, we see in \Cref{fig:MaterialSolutions} that the largest global errors in the gradient are limited to the boundary, and the error is dispersed uniformly over the domain.

We now aim to make a comparison with a classical feed-forward fully-connected architecture. As the jumps in the gradient cannot be appropriately defined for smooth architectures, we cannot use a PINN-based loss to describe the system. {\blue Thus,} in order to make a fair comparison between the two architectures, we employ a discretisation of the $H^1$-norm as a loss function, as we described in Section \ref{secMotivation}.

For the regularity-conforming architecture, we consider the same architecture as in the {\blue previous} PINNs implementation, {\blue yielding 2710 trainable parameters}. For the classical architecture, we employ a fully-connected feed-forward NN with $tanh$ activation function and 3 hidden layers of 36 neurons each, {\blue giving a total of 2809 trainable parameters.} 

\begin{figure}[H]\begin{center}
\begin{subfigure}[t]{0.45\textwidth}\begin{center}
\includegraphics[height=0.85\textwidth]{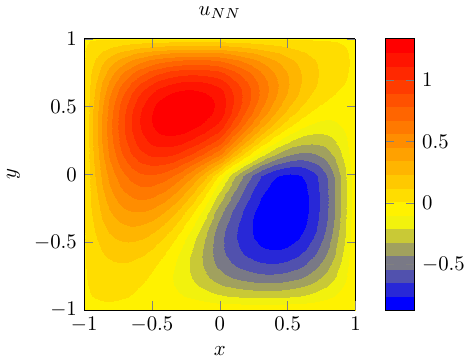}
\caption{Approximate solution .}
\end{center}\end{subfigure}
\begin{subfigure}[t]{0.45\textwidth}\begin{center}
\includegraphics[height=0.85\textwidth]{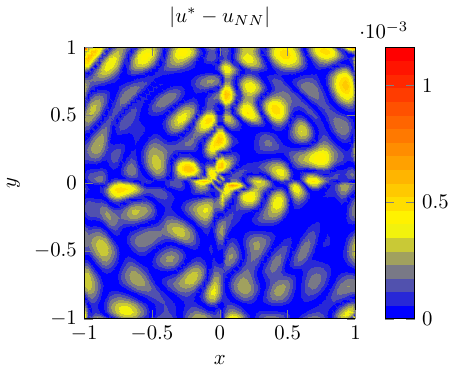}
\caption{Pointwise error .}
\end{center}\end{subfigure}
\begin{subfigure}[t]{0.45\textwidth}\begin{center}
\includegraphics[height=0.85\textwidth]{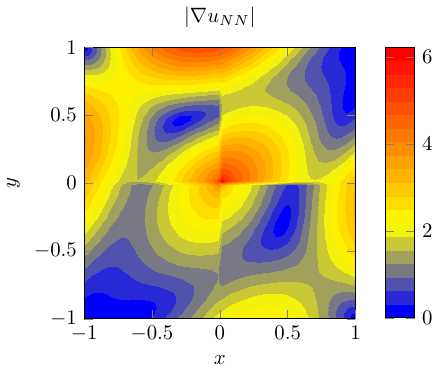}
\caption{Norm of the gradient of the approximate solution }
\end{center}\end{subfigure}
\begin{subfigure}[t]{0.45\textwidth}\begin{center}
\includegraphics[height=0.85\textwidth]{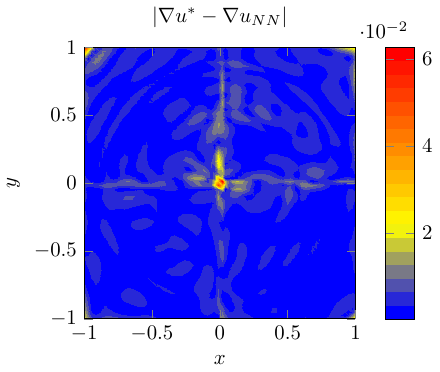}
\caption{Error of the gradient of the approximate solution }
\end{center}\end{subfigure}\caption{The approximate solutions for a conforming architecture trained on the $H^1$-norm and errors. $u_{NN}$ has a final relative error of $0.033\%$ in $L^2$ and an $L^2$ error in the gradient of $0.24\%$. }\label{fig:MaterialH1SolutionsConforming}
\end{center}
\end{figure}

\begin{figure}[H]\begin{center}
\begin{subfigure}[t]{0.45\textwidth}\begin{center}
\includegraphics[height=0.85\textwidth]{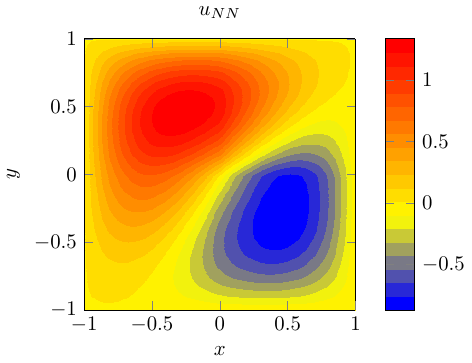}
\caption{Approximate solution $u_{NN}$ \\(Classical architecture).}
\end{center}\end{subfigure}
\begin{subfigure}[t]{0.45\textwidth}\begin{center}
\includegraphics[height=0.85\textwidth]{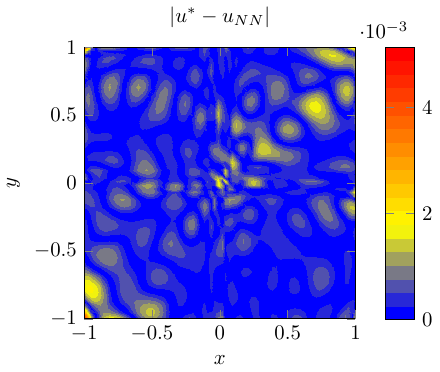}
\caption{Pointwise error \\(Classical architecture).}
\end{center}\end{subfigure}
\begin{subfigure}[t]{0.45\textwidth}\begin{center}
\includegraphics[height=0.85\textwidth]{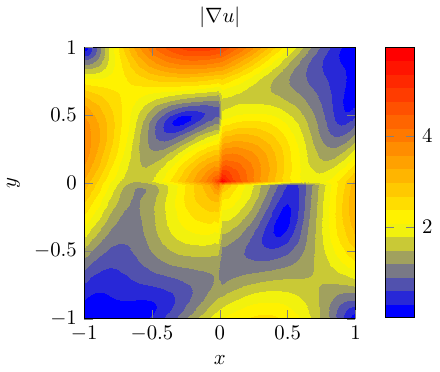}
\caption{Norm of the gradient of the approximate solution \\ (Classical architecture)}
\end{center}\end{subfigure}
\begin{subfigure}[t]{0.45\textwidth}\begin{center}
\includegraphics[height=0.85\textwidth]{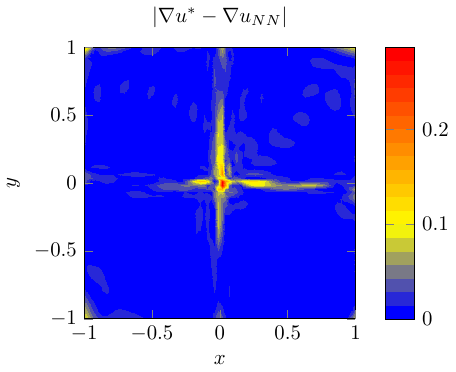}
\caption{Error of the gradient of the approximate solution \\(Classical architecture)}
\end{center}\end{subfigure}
\caption{{The approximate solutions for a classical architecture trained on the $H^1$-norm and errors. $u_{NN}$ has a final relative error of $0.08\%$ in $L^2$ and a relative $L^2$-error in the gradient of $1.57\%$}}
\label{fig:MaterialH1SolutionsClassical}
\end{center}
\end{figure}

\begin{figure}[H]\begin{center}
\includegraphics[width=0.5\textwidth]{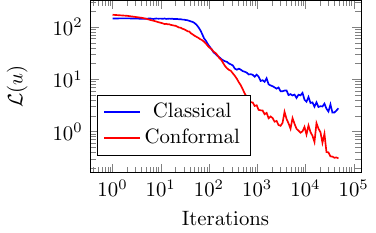}
\caption{Loss evolution during training for the Classical and Conforming architectures.}\label{fig:MaterialH1Loss}\end{center}
\end{figure}

We employ as a loss function in both cases 
\begin{equation}
\mathcal{L}(u_{NN})=\left(\frac{1}{N}\sum\limits_{i=1}^N |u(x_i)-u_{NN}(x_i)|^2+|\nabla u(x_i)-\nabla u_{NN}(x_i)|^2\right)^\frac{1}{2}, 
\end{equation}
where $\{x_i\}_{i=1}^N$ are uniformly sampled {\blue points} from $\Omega$ and randomly generated at each iteration. For the implementation, we consider $N=25,00$, and use the Adam optimiser for 50,000 iterations, using a learning rate of $10^{-3}$ for the first 25,000 iterations and exponentially decaying to $10^{-6}$ for the last 25,000 iterations. Figures \ref{fig:MaterialH1SolutionsConforming} and \ref{fig:MaterialH1SolutionsClassical} show the approximation results of the conforming and classical architectures, respectively, and Figure \ref{fig:MaterialH1Loss} shows the evolution of the loss during training.


We observe a good $L^2$ approximation in both cases, which is understandable as the target function is continuous and the loss explicitly controls the $L^2$ error. With the classical architecture, we obtain a relative $L^2$-error of 0.08\% and a relative $L^2$-error in the gradient of 1.57\%. Whereas for the conforming architecture, the relative $L^2$-error is 0.033\% and the relative $L^2$-error in the gradient of 0.24\%. In the case of the conforming architectures, there are localised errors in the vicinity of the jumps in the gradient and a localised error at the singularity. Nonetheless, these are significantly milder than those obtained using the classical architecture, being of order $10^{-2}$ in the case of the regularity-conforming architecture and an order of magnitude higher with the classical architecture. We remark that in each case, the error in the gradient at the origin is unbounded, and the maximal values ($0.06$ and $0.25$, respectively) are representative of the finite resolution of the graphs. 

\section{Conclusions}\label{secConclusions}
This work introduces Regularity-Conforming Neural Networks (ReCoNNs) to approximate solutions of 2D transmition problems presenting different types of singularities. We define such architectures by employing \textit{a priori} information on the location and nature of the singular behaviour, while the exact form of the singularity is approximated during the training process.
We consider loss functions in strong form similar to the ones considered in PINNs. These architectures allow us to easily incorporate interface conditions in the loss functional as they enjoy partial explainability. We test our method in problems exhibiting gradient discontinuities and point singularities coming from re-entrant corners and vertices {\blue at} material interfaces. We show the superiority of employing ReCoNNs in contrast to classical feed-forward NNs as they overcome the issue of Gibbs-type phenomena and instabilities near the singularity, leading to a good approximation and stable convergence towards the solution.

Possible future research line include: (a) the application of ReCoNNs to 3D problems (i.e. the Fichera {\blue problem}), (b) the generalization of the method to more challenging problems like crack singularities or multiply-connected domains, and (c) the use of ReCoNNs infrastructure to solve inverse problems to characterize the material properties involved in transmission problems.

\section*{Acknowledgements}
David Pardo has received funding from: the Spanish Ministry of Science and Innovation projects with references TED2021-132783B-I00 funded by MCIN/ AEI /10.13039/501100011033 and the European Union NextGenerationEU/ PRTR, PID2019-108111RB-I00 funded by MCIN/ AEI /10.13039/501100011033, the “BCAM Severo Ochoa” accreditation of excellence CEX2021-001142-S / MICIN / AEI / 10.13039/501100011033; the Spanish Ministry of Economic and Digital Transformation with Misiones Project IA4TES (MIA.2021.M04.008 / NextGenerationEU PRTR); and the Basque Government through the BERC 2022-2025 program, the Elkartek project BEREZ-IA (KK-2023/00012), and the Consolidated Research Group MATHMODE (IT1456-22) given by the Department of Education. Judit Muñoz-Matute has received funding from the European Union’s Horizon 2020 research and innovation programme under the Marie Sklodowska-Curie individual fellowship grant agreement No. 101017984 (GEODPG).
\bibliography{bib}

\appendix

\section{Construction of the exact solution for the interior vertex problem}\label{App}

In Section \ref{subsecVertices} we considered the transmission problem on $(-1,1)^2$ with 
\begin{equation}
\sigma(x)=\left\{\begin{array}{c c}
1 & x_1>0,x_2>0,\\
2 & x_1<0,x_2>0,\\
3 & x_1<0,x_2<0,\\
4 & x_1>0,x_2<0. 
\end{array}\right.
\end{equation}
Herein, we outline the process we used to construct $u^*$, the exact solution to the problem. 

First, we consider the Sturm-Liouville equation, associated with the singular function given in polar coordinates. With abuse of notation, we identify $\sigma$ with its representation in polar coordinates centred at $0$. The Sturm-Liouville equation satisfied by the singular function is thus 
$$
\int_0^{2\pi}\sigma(\theta)\left(\phi'(\theta)v'(\theta)-\lambda^2\phi(\theta)v(\theta)\right)\,d\theta=0
$$
for all $v\in H^1_{per}(0,2\pi)$, the space of all $2\pi$-periodic functions in $H^1$. It is then immediate to see that any $\phi$ satisfying this equation may be expressed piecewise as 
\begin{equation}
\phi(\theta)=\left\{\begin{array}{c c}
a_1\sin(\lambda \theta)+b_1\cos(\lambda\theta) & 0<\theta<\frac{\pi}{2},\\
a_2\sin(\lambda \theta)+b_2\cos(\lambda\theta) & \frac{\pi}{2}<\theta<\pi,\\
a_3\sin(\lambda \theta)+b_3\cos(\lambda\theta) & \pi<\theta<\frac{3\pi}{2},\\
a_4\sin(\lambda \theta)+b_4\cos(\lambda\theta) & \frac{3\pi}{2}<\theta<2\pi. 
\end{array}\right.
\end{equation}
There are 9 unknowns, corresponding to the coefficients $a_i,b_i$ and $\lambda$. As $\phi\in H^1_{per}$, $\phi$ is continuous, and we thus obtain four linear equations in the coefficients $a_i,b_i$, each of the form 
$$
\lim\limits_{\theta\to\theta_i^-}\phi(\theta)=\lim\limits_{\theta\to\theta_i^+}\phi(\theta),
$$
where $\theta_i$ is either $0,\frac{\pi}{2},\pi$ or $\frac{3\pi}{4}$. Similarly, the continuity of the flux gives a further four linear equations in the coefficients $a_i,b_i$ via 
$$
\lim\limits_{\theta\to\theta_i^-}\sigma(\theta)\phi'(\theta)=\lim\limits_{\theta\to\theta_i^+}\sigma(\theta)\phi'(\theta).
$$
This yields a total of 8 equations in 9 unknowns. For a fixed but arbitrary $\lambda$, we may then represent the equations as an $8\times 8$ matrix $M_\lambda$, depending non-linearly on $\lambda$, acting on the coefficients. The existence of non-trivial solutions to the Sturm-Liouville equation thus requires that $\det(M_\lambda)=0$. Any standard root-finding algorithm may be applied, and we obtain a single root in $(0,1)$. Once this root is obtained, it suffices to find the kernel of $M_\lambda$ at this value in order to obtain the singular solution. 

Once the singular solution, which we denote $s^*(x)$, has been obtained, we then define $u^*(x)=\cos\left(\frac{x_1\pi}{2}\right)\cos\left(\frac{x_2\pi}{2}\right)s^*(x)$. We then define our forcing term pointwise $f(x)=\sigma(x)\Delta u^*(x)$ and verify that $u^*$ satisfies the weak-form equation 
$$
\int_\Omega \sigma(x)\nabla u^*(x)\cdot \nabla v(x)+f(x)v(x)\,dx=0
$$
for all $v\in H^1_0(\Omega)$. From the construction, it is clear that it is only necessary to verify the continuity of the flux at the interfaces. Considering the interface described by $x_1=0$, we have that the normal derivative corresponds to the derivative with respecto to $x_1$, and thus 
\begin{equation}
\begin{split}
\frac{\partial u^*}{\partial \nu}(x_1,x_2)=\frac{\pi}{2}\sin\left(\frac{x_1\pi}{2}\right)\cos\left(\frac{x_2\pi}{2}\right)s^*(x)+\cos\left(\frac{x_1\pi}{2}\right)\cos\left(\frac{x_2\pi}{2}\right)\frac{\partial s^*}{\partial \nu}(x).
\end{split}
\end{equation}
At the interface $x_1=0$, the first term vanishes, and thus the continuity of the flux in $u^*$ is a direct consequence of the continuity of the flux in $s^*$, as 
\begin{equation}
\sigma(x)\frac{\partial u^*}{\partial\nu}=\cos\left(\frac{x_1\pi}{2}\right)\cos\left(\frac{x_2\pi}{2}\right)\left(\sigma(x)\frac{\partial s^*}{\partial \nu}(x)\right).
\end{equation}
The same argument may then be applied to the interface corresponding to $x_2=0$, and we conclude that $u^*$ satisfies the corresponding weak-form transmission problem. 
\end{document}